\documentclass[12pt]{article}

\newif\iffrench\frenchfalse
\usepackage[utf8]{inputenc}
\usepackage[T1]{fontenc}

\iffrench
\usepackage[french]{babel}
\else
\usepackage[english]{babel}
\fi

\usepackage{pslatex}
\usepackage{lmodern}

\usepackage{amssymb}
\usepackage{amsthm}
\usepackage{bbold} % indicatrices
\usepackage{dsfont}
\usepackage{mathtools}
\usepackage{thmtools}

\usepackage{varioref} % Intelligent page references
\usepackage{hyperref} % Extensive support for hypertext in LaTeX
\hypersetup{
    colorlinks=true,
    linkcolor=blue,
    }

\usepackage{tikz}
\usepackage{pstricks-add}

\usepackage{enumitem} % Control layout of itemize, enumerate, description
\usepackage[babel=true]{csquotes} % Context sensitive quotation facilities
\usepackage{color}
\usepackage{setspace}
\usepackage{array}
\usepackage{layout}

\usepackage{stmaryrd} % St Mary Road symbols for theoretical computer science
\usepackage{mathrsfs} % jolies cursives (mathscr)

\usepackage[margin=20pt,font=small,labelfont=bf]{caption} % jolis titres de figures

\usepackage[
  backend=biber,
  citestyle=numeric-comp,
  giveninits=true,
  url=false,
  isbn=false,
  doi=false,
  eprint=false,
  maxbibnames=50,
]{biblatex}
\newcommand{\smarge}[2]{\usepackage[top=#1,bottom=#1+1cm,left=#1-#2,right=#1]{geometry}}

\iffrench

\addto\captionsfrench{} % nom des tableaux en français

\newcommand{\ioo}[1]{\left]#1\right[}
\newcommand{\ifo}[1]{\left[#1\right[}
\newcommand{\iof}[1]{\left]#1\right]}

\newtheorem{thm}{Theorème}[section]
\newtheorem{ppn}[thm]{Proposition}
\newtheorem{cor}[thm]{Corollaire}
\newtheorem{lem}[thm]{Lemme}
\newtheorem{dfi}[thm]{Définition}
\newtheorem{cjt}[thm]{Conjecture}
\newcommand{\pth}[1]{\left(#1\right)}
\newcommand{\cro}[1]{\left[#1\right]}
\newcommand{\acc}[1]{\left\{#1\right\}}
\newcommand{\abs}[1]{\left|#1\right|}
\newcommand{\dabs}[1]{\left\|#1\right\|}

\newcommand{\ceil}[1]{\left\lceil#1\right\rceil}

\else

\newtheorem{thm}{Theorem}[section]
\newtheorem{ppn}[thm]{Proposition}

\newtheorem{lem}[thm]{Lemma}
\newtheorem{dfi}[thm]{Definition}

\theoremstyle{remark}
\newtheorem{req}[thm]{Remark}

\newcommand{\ioo}[1]{\left(#1\right)}
\newcommand{\ifo}[1]{\left[#1\right)}
\newcommand{\iof}[1]{\left(#1\right]}

\newcommand{\pth}[1]{\left(#1\right)}
\newcommand{\cro}[1]{\left[#1\right]}
\newcommand{\acc}[1]{\left\{#1\right\}}
\newcommand{\abs}[1]{\left|#1\right|}
\newcommand{\dabs}[1]{\left\|#1\right\|}

\newcommand{\ceil}[1]{\left\lceil#1\right\rceil}

\newcommand{\ie}{i.e. } 
\newcommand{\eg}{e.g. }

\fi
\newcommand{\esp}{\hspace{1cm}}
\newcommand{\comment}[1]{\hspace{0.5cm}\text{#1}\hspace{0cm}}

\newcommand{\tq}{\hspace{0.25cm}/ \hspace{0.25cm}}
\newcommand{\vg}{,\,}

\newcommand{\goq}{\geqslant}
\newcommand{\loq}{\leqslant}
\newcommand{\eps}{\varepsilon}

\newcommand{\sys}[1]{\begin{equation}\left\{\begin{aligned}#1\end{aligned}\right.\end{equation}}

\newcommand{\pmat}[1]{\begin{pmatrix} #1\end{pmatrix}}

\newcommand{\de}{\,\mathrm{d}}

\newcommand{\dr}{\partial}

\newcommand{\toweak}{\rightharpoonup}

\newcommand{\Er}{\mathds{R}}

\newcommand{\mcc}{\mathcal{C}}

\newcommand{\mce}{\mathcal{E}}

\newcommand{\mcl}{\mathcal{L}}

\newcommand{\mct}{\mathcal{T}}

\smarge{2.5cm}{0cm}

\renewbibmacro{in:}{}

\title{Qualitative properties of the spreading speed of a population structured in space and in phenotype}

\author{Nathanaël Boutillon \\
\footnotesize{Aix Marseille Université, CNRS, I2M, Marseille, France}\\
\footnotesize{INRAE, BioSP, 84914, Avignon, France}\\
}

\begin{document}

\maketitle

\begin{abstract}
  We consider a nonlocal Fisher-KPP equation that models a population structured in space and in phenotype. The population lives in a heterogeneous periodic environment: the diffusion coefficient, the mutation coefficient and the fitness of an individual may depend on its spatial position and on its phenotype.

  We first prove a Freidlin-Gärtner formula for the spreading speed of the population. We then study the behaviour of the spreading speed in different scaling limits  (small and large period, small and large mutation coefficient). Finally, we exhibit new phenomena arising thanks to the phenotypic dimension. Our results are also valid when the phenotype is seen as another spatial variable along which the population does not spread.
\end{abstract}



\tableofcontents

\vspace{2cm}

\section{Introduction}\label{s:intro}

Our goal is to study the spreading properties of a population structured in space and in phenotype, and which evolves according to a reaction-diffusion equation.
A particular case of our model reads: 
\begin{equation}\label{eq:main_intro}
  \left\{
  \begin{aligned}
    \dr_tu&=a\Delta_x u+m\Delta_{\theta}u+u\pth{r(x,\theta)-\rho(t,x)},&t>0\vg(x,\theta)\in\Er\times\Theta,\\
    \nu\cdot\nabla_{\theta} u&=0,&t>0\vg(x,\theta)\in\Er\times\dr\Theta,\\
    u(0,x,\theta)&=u_0(x,\theta),&(x,\theta)\in\Er\times\Theta,
  \end{aligned}
  \right.
\end{equation}
with
\[\rho(t,x):=\int_{\Theta}u(t,x,\sigma)\de \sigma.\]
Here~$\Theta\subset\Er^P$ is a smooth open bounded domain with outward normal~$\nu$.
The initial condition~$u_0$ is assumed to be nonnegative, continuous and bounded.

We call~$x\in\Er$ the space variable and~$\theta\in\Theta$ the phenotype variable,
and we see $u(t,x,\theta)$ as the density of individuals of phenotype~$\theta\in\Theta$ living at position~$x\in\Er$ and at time~$t\goq 0$.
%If~$u_0\goq 0$, then~$u\goq 0$ for all~$t>0$ and we can see~$u(t,x,\theta)$ as the density of individuals of phenotype~$\theta$ living at position~$x\in\Er$ and at time~$t\goq 0$.
With this viewpoint,~$r(x,\theta)$ is the fitness of an individual of phenotype~$\theta$ at position~$x$, while $a>0$ is the motility and~$m>0$ is the mutation coefficient ($m$ corresponds to the product between the mutation rate and the mutational variance, see \eg~\cite[Appendix A]{HLMR20}). Later, we shall rather see~$m$ as a scaling parameter.
We assume that the environment is periodic in space, which we translate into our framework by assuming that~$r$ is periodic in the variable~$x$.

We will show that a population modelled by Equation \eqref{eq:main_intro} invades its environment at some finite speed, which we will call the \emph{spreading speed}. We will compute the spreading speed for different scaling limits: when the period of the environment is large or small, and when the mutation coefficient is large or small. Finally, we will show that there exists a mutation coefficient~$m$ which maximises the spreading speed. %From the point of view of modelling, both problems are important.

\paragraph{Simpler models.}
Before describing our precise assumptions and our results, let us explain where Equation~\eqref{eq:main_intro} comes from.
The Fisher-KPP equation was introduced in 1937 by Fisher~\cite{Fis37} and Kolmogorov, Petrovsky and Piskounov~\cite{KPP37}. It reads:
\begin{equation}\label{eq:kpp_classique}
  \dr_tu(t,x)=a\Delta_x u+ru-u^2,
\end{equation}
where~$a>0$ and~$r>0$ are constant.
A motivation for Equation \eqref{eq:kpp_classique} comes from biology:~$u(t,x)$ can be seen as the density of individuals living at position~$x\in\Er$ and at time~$t\goq 0$. The term~$a\Delta_x u$ means that the individuals move according to a (time-changed) Brownian motion. The linear term~$ru$ stands for the demography in the absence of interactions. The population is regulated logistically through the nonlinear term~$-u^2$, which accounts for the additional deaths due to the competition between individuals. 

One very important feature of the Fisher-KPP equation is the emergence of invasions of the environment by the population. Namely, there exists a \emph{spreading speed}~$c_{KPP}:=2\sqrt{ar}$ satisfying the following property. If the initial condition is nonnegative, is nonzero and has a support bounded from above, then, for all~$c'>c_{KPP}$,
\[\lim_{t\to+\infty}u(t,c't)=0,\]
and for all~$c''\in\ioo{0,c_{KPP}}$,
\[\lim_{t\to+\infty}u(t,c''t)=r.\]
Moreover, for all~$c\goq c_{KPP}$, there exists a profile~$U_c$ such that~$u(t,x):=U_c(x-ct)$ solves \eqref{eq:kpp_classique}. Such a self-similar solution is called a \emph{travelling wave}.

The simplicity of the Fisher-KPP equation, together with these invasion properties, have made it a central tool in ecology and in biology (\eg~\cite{Ske51,Tur98,ShiKaw97}).
However, because of its simplicity, the Fisher-KPP equation misses phenomena which may play an important role in invasions. Equation \eqref{eq:main_intro} is a variant that tries to take two such phenomena into account: the heterogeneity of the environment and the heterogeneity of the population.

\paragraph{Heterogeneous environment.}
A first phenomenon which is missed by the classical Fisher-KPP equation is the heterogeneity of the environment. This can be (partially) overcome by considering the following periodic heterogeneous version of the Fisher-KPP equation:
\begin{equation}\label{eq:kpp_hetero}
  \dr_tu(t,x)=a\Delta_x u+r(x)u-u^2,
\end{equation}
where the function~$r$ is smooth and~$1$-periodic.
In this setting, if the population persists, there still exists a minimal spreading speed~$c_{het}$ (which plays a role analogous to~$c_{KPP}$ but in the heterogeneous setting).
Gärtner and Freidlin~\cite{GarFre79} found a variational formula for~$c_{het}$.  
For~$\lambda>0$, they consider the operator
\[\mcl^{\lambda}:=a\Delta_x-2\lambda a\dr_x+(r(x)+\lambda^2a)\]
acting on~$1$-periodic functions.
The principal eigenvalue~$k^{\lambda}$ of this operator is the unique real number such that there exists a~$1$-periodic function~$\phi>0$ satisfying ~$\mcl^{\lambda}\phi=k^{\lambda}\phi$. 
Gärtner and Freidlin show that if~$k^0>0$, then the population persists and an invasion occurs at speed
\begin{equation}\label{eq:freidlin_gartner_intro}
  c_{het}:=\inf_{\lambda>0}\frac{k^{\lambda}}{\lambda}.
\end{equation}
Since Gärtner and Freidlin, numerous work have proved different versions of Equation \eqref{eq:freidlin_gartner_intro} for various settings, \eg~\cite{BHN05,BHN08,Wei02,Gir18,GriMat21}.
In the heterogeneous setting, the existence of travelling waves cannot be expected, due to the heterogeneity. Still, there exist \emph{pulsating travelling waves} for all speeds no less than the minimal spreading speed~$c_{het}$. A pulsating travelling wave with speed~$c\goq c_{het}$ is a solution~$v$ of \eqref{eq:kpp_hetero} such that~$(t,x)\mapsto v(t,x+ct)$ is periodic in time and which has the following limits: for all~$c'>c_{het}$,~$v(t,c't)$ converges to $0$, and for all~$c''<c_{het}$, the profile~$x\mapsto v(t,c''t+x)$ converges locally, up to a time-dependent shift, to the unique~$\phi>0$ satisfying~$a\Delta_x\phi+r(x)\phi-\phi^2=0$. The limiting profile~$\phi$ is called the stationary state of the equation and, by uniqueness, is~$1$-periodic.
See also~\cite{BerHam02,BHR05-2}. In~\cite{Xin00}, Xin gives a survey of results and methods for reaction-diffusion equations in heterogeneous environments.

\paragraph{Heterogeneous population.}  A second phenomenon which is missed by the basic Fisher-KPP equation \eqref{eq:kpp_classique} is the heterogeneity of the population, \ie some individuals may react differently to the same environmental conditions. Thus we add a new variable~$\theta\in\Theta$ representing the phenotype of the individuals. If we do this directly on \eqref{eq:kpp_hetero}, we get
\begin{equation}\label{eq:kpp_hetero2}
  \dr_tu(t,x,\theta)=a\Delta_xu+m\Delta_{\theta}u+r(x,\theta)u-u^2,\esp x\in\Er\vg \theta\in \Theta,
\end{equation}
together with the boundary condition~$\nu\cdot\nabla_{\theta} u=0$ on~$\dr \Theta$. Here~$m\Delta_{\theta}$ accounts for the mutations, with a mutation coefficient~$m$.
Then the Freidlin-Gärtner formula also holds~\cite{Wei02}, and pulsating travelling waves also exist~\cite{BHN05,BHN08}. 

The only difference between Equation \eqref{eq:main_intro} and Equation \eqref{eq:kpp_hetero2} is the nonlocal term~$\rho$. The nonlocal term~$\rho$ comes from the fact that in our model, the competition to which an individual is subject should depend on \emph{all} the individuals that are located at the same position, regardless of their phenotypes. Thus the competition term~$-u^2$ must be replaced by~$-\rho u$.
An equation similar to \eqref{eq:main_intro} (with rarer but larger mutations) appeared for the first time in the work of Prévost~\cite{Pre04} and was then derived by Champagnat and Méléard~\cite{ChaMel07} as a large population limit of an individual-based model.

The discussion above implies that \eqref{eq:kpp_hetero2} is more adapted to populations that depend on variables~$x$ and~$\theta$ which are of the same nature (for example,~$x$ and~$\theta$ could be the coordinates of a spatial position), while we should better use \eqref{eq:main_intro} if we want the variable~$\theta$ to represent a phenotype. Although we shall concentrate on Equation \eqref{eq:main_intro}, our results only depend on the linear problem, so they are valid for the local equation~\eqref{eq:kpp_hetero2} as well (in the latter setting, however, our results may be harder to interpret from a biological viewpoint).

For the sake of exposition, let us assume temporarily that the fitness~$r$ only depends on the phenotype~$\theta$ in~\eqref{eq:main_intro}, so that the environment becomes homogeneous (\ie independent of the spatial position). Then, there may be individuals, say with phenotype~$\theta_1$, which are more favoured than others, say with phenotype~$\theta_2$. It can then be expected that individuals of phenotype~$\theta_1$ will drive the invasions, while individuals of phenotype~$\theta_2$ will slow it down. Let us now try to understand the influence of this heterogeneity on the spreading speed.
Let~$H$ be the principal eigenvalue of the operator
\[m\Delta_{\theta}+r(\theta)\]
acting on the functions~$\psi\in\mcc^2(\overline{\Theta})$ such that~$\nu\cdot\nabla_{\theta}\psi=0$ on~$\dr\Theta$. Namely,~$H$ is the only real number such that there exists a function~$\psi>0$ satisfying
\[m\Delta_{\theta}\psi+r(\theta)\psi=H\psi,\esp\nu\cdot\nabla_{\theta}\psi=0.\]
Assume~$H>0$. Then, if~$u_0$ is nonnegative, not uniformly equal to~$0$ and has a support bounded from above, then the solution of \eqref{eq:main_intro} spreads at a speed~$c_{phenotype}=2\sqrt{aH}$, in the sense that for all~$c'>c_{phenotype}$,
\[\lim_{t\to+\infty}\rho(t,c't)=0,\]
and for all~$c''\in\ioo{0,c_{phenotype}}$,
\[\liminf_{t\to+\infty}\rho(t,c''t)>0.\]
Moreover, for all~$c\goq c_{phenotype}$, there exists a profile~$U_{c}$ such that~$u(t,x,\theta):=U(x-ct,\theta)$ solves \eqref{eq:main_intro}. See~\cite{BJS14}. Note that $H$ plays a role similar to $r$ for the usual Fisher-KPP speed $c_{KPP}$.

We now go back to the general case of a heterogeneous environment and a heterogeneous population with a continuous phenotype. %, and survey the recent literature about it.

\paragraph{Heterogeneous environment and population.} 
In~\cite{ACR13,Pel20}, the authors studied the case where the population evolves along an environmental gradient, that is:
\[r(x,\theta)=R(x-B\theta),\esp B>0,\]
with
\begin{align*}
  R>0&\comment{on $\ioo{-1,1}$},\\
  R<0&\comment{on $\Er\setminus\cro{-1,1}$}.
\end{align*}
In this model, the phenotype of the individuals must be higher and higher as the population moves to the right: evolution is necessary for the propagation of the population. This is relevant when we want to study a propagation along the north-south axis, or uphill or downhill if the population lives on a mountain. In~\cite{Pel20}, Peltier showed that if the tail of the initial condition is heavy enough, then an accelerating invasion arises, as in the~$1$-dimensional case~\cite{HamRoq10}. In~\cite{ACR13}, Alfaro, Coville and Raoul gave a criterion, based on the principal eigenvalue of an elliptic operator, to determine whether the population persists or gets extinct. They also showed the existence of a travelling wave.

The authors of~\cite{ABR17} added a heterogeneity in time corresponding to climate change. In this new setting, the growth rate~$r(x,\theta)$ is replaced by~$r(t,x,\theta)=R(x-ct,\theta)$.
They study different shapes of the favourable space-phenotype zone (confined zone, environmental gradient, or a mixing between a confined zone and an environmental gradient). In each case, they found a criterion for the persistence of the species, and found a formula for the spreading speed.

An interesting case of equations with phenotype is the \emph{cane-toad equation}, where~$\Theta\subset\ioo{0,+\infty}$ may be unbounded, and the diffusivity~$a$ is nonconstant and proportional to the phenotype:~$a(\theta)=\theta$. In this case, accelerating invasions may arise~\cite{BCMMPRV12}. Then the interplay with the phenotype has a striking behaviour with respect to the spreading properties and real world biological consequences. In this work, however, most results do not apply to the cane-toad equation.

There are also works about periodic environments, as modelled by our equation~\eqref{eq:main_intro}. In~\cite{AlfPel22}, Alfaro and Peltier study the effect of small perturbations of the homogeneous case on the spreading speed and the travelling waves. They also deal with a non-gradient and non-periodic heterogeneity.
In~\cite{BouMir15}, Bouin and Mirrahimi study the effect of large periods and get the precise position of the invasion, using Hamilton-Jacobi equation techniques.
  In~\cite{LecMir23}, Léculier and Mirrahimi study the weak mutation limit of the equilibria of an equation which is related to ours, but displays nonlocal mutations.
In~\cite{br24}, Boutillon and Rossi focus on the same equation~\eqref{eq:main_intro} but possibly set in an unbounded phenotype space. They prove the well-posedness of~\eqref{eq:main_intro} and give a criterion, based on a generalised principal eigenvalue, for persistence.

As a concluding remark, we mention that another way to study heterogeneous populations is to use a system of Fisher-KPP equations: each equation of the system represents a variant which competes with the other ones. We mention~\cite{Gir18} for homogeneous environments, and~\cite{Gir23,GriMat21,rbzp24} for heterogeneous environments. In the latter, Freidlin-Gärtner formulas similar to ours (but in discrete-phenotype settings) are proved and their properties are studied.

\paragraph{Problematic and layout.} Our main goal is to study the spreading speed of a population modelled by equations of the form \eqref{eq:main_intro}. In the next section, we present our results. We first give our general assumptions and the general form of the equation on which we shall work. In Subsection~\ref{ss:fg}, we shall give a variational formula analogous to the Freidlin-Gärtner formula. 

From a biological point of view, it is important to note that movements of individuals and mutations are phenomena of different nature; therefore, there is no reason why they should occur on the same timescales. It is thus important to understand the interplay between the mutation coefficient and the period at different scales. What is the spreading speed for large and small periods? What is the spreading speed for large and small mutation coefficients? Subsection~\ref{ss:different_scales} is devoted to these questions.  We shall see that in these limits, one can often get back to the study of equations without the phenotypic dimension.

A specificity of Equation \eqref{eq:main_intro} is that the population modelled by the equation is heterogeneous. What are the new phenomena which can occur? What is their biological interpretation? In Subsection~\ref{ss:mutation_effect}, we shall give several results which are invisible when we forget the phenotypic dimension.

\section{Main results}\label{s:main_results}

Now, we describe our general model and our general assumptions.
We consider a mutation coefficient~$m>0$ and a period~$L>0$. 
We let~$\Theta$ be an open bounded domain with~$\mcc^{2}$ boundary.

We fix~$\alpha\in\ioo{0,1}$.
We take $a\in\mcc^{1,\alpha}(\Er\times\overline{\Theta})$, $\mu\in\mcc^{1,\alpha}(\Er\times\overline{\Theta})$ and $r\in\mcc^{0,\alpha}(\Er\times\overline{\Theta})$ such that for all~$\theta\in\overline{\Theta}$,
the functions $x\mapsto a(x,\theta)$, $x\mapsto \mu(x,\theta)$ and $x\mapsto r(x,\theta)$ are~$1$-periodic (we say that $a$, $\mu$ and $r$ are \emph{$1$-periodic in $x$}).

The functions $a$ and $\mu$ represent the heterogeneity of the diffusion in the directions $x$ and $\theta$ respectively.
  We make an assumption of uniform ellipticity of the operator, that is,
we assume that there exists a constant $\eta>0$ such that $a\goq\eta$ and $\mu\goq\eta$.

As we shall be interested in the behaviour of the spreading speed for various scales of the period and of the mutation coefficient, we first write out our model in a way that makes explicit the dependence in these parameters. We note
\[a_L(x,\theta):=a\pth{\frac{x}{L},\theta},\esp \mu_L(x,\theta):=\mu\pth{\frac{x}{L},\theta},\esp r_L(x,\theta):=r\pth{\frac{x}{L},\theta}\]
the~$L$-periodic versions of~$a$,~$\mu$ and~$r$ respectively. We define an operator~$\mct_{m,L}$ by
\begin{equation}\label{eq:def_mct}
  \mct_{m,L}\phi:=\dr_x\pth{a_L(x,\theta)\dr_x\phi}+m\nabla_{\theta}\cdot\pth{\mu_L(x,\theta)\nabla_{\theta}\phi}+r_L(x,\theta)\phi,
\end{equation}
and we write a general version of main model \eqref{eq:main_intro} as:
\begin{equation}\label{eq:main}
  \left\{
  \begin{aligned}
    \dr_tu_{m,L}&=\mct_{m,L}u_{m,L}-u_{m,L}\times\rho_{m,L}(t,x),&t>0\vg(x,\theta)\in\Er\times\Theta,\\
    \nu\cdot\nabla_{\theta} u_{m,L}&=0,&t>0\vg (x,\theta)\in\Er\times\dr\Theta,\\
    u_{m,L}(0,x,\theta)&=u_0(x,\theta),&(x,\theta)\in\Er\times\Theta,
  \end{aligned}
  \right.
\end{equation}
with
\[\rho_{m,L}(t,x):=\int_{\Theta}u_{m,L}(t,x,\sigma)\de\sigma.\]
We also require that $u_0\in\mcc^{0}(\Er\times\overline{\Theta})$ is nonnegative, bounded and has a support bounded from above.
Contrarily to \eqref{eq:main_intro}, $a$ and~$\mu$ are now allowed to depend on~$x$ and~$\theta$.
These dependencies will be restricted in several of our results below.
%We see the variable~$m$ as a scaling parameter.
%Since the heterogeneity $\mu$ is fixed, we will for convenience call $m$ the mutation coefficient (instead of $m\mu$). 

Theorem~2.1 of~\cite{br24} ensures that the problem \eqref{eq:main} has a unique weak solution in $W^{1,2}_{p,loc}(\ifo{0,+\infty}\times(\Er\times\overline{\Theta}))$ for all $p>1$
(\ie the solution, its first derivatives with respect to $t$, and all its derivatives with respect to $(x,\theta)$ up to the second order, belong to $L^p_{loc}$).
Moreover, the unique solution $u$ and the corresponding $\rho$ are continuous, globally bounded and nonnegative. If $u_0$ is nonzero, then $u$ and $\rho$ are strictly positive at all positive times.

\subsection{The Freidlin-Gärtner formula}\label{ss:fg}

We first introduce a family of principal eigenvalues which will be central in the statement of the Freidlin-Gärtner formula. For $\lambda>0$, we let~$\mcl^{\lambda}_{m,L}$ be the linear operator defined by:
\[\mcl^{\lambda}_{m,L}:=\mct_{m,L}-2\lambda a_L\dr_x+\pth{\lambda^2a_L-\lambda\dr_xa_L},\]
 or, equivalently,
 \[\mcl^{\lambda}_{m,L}\phi(x,\theta):=e^{\lambda x}\mct_{m,L}(e^{-\lambda \cdot}\phi(\cdot,\theta)).\]
 We let $k^{\lambda}_{m,L}\in\Er$ be the principal eigenvalue associated to the eigenproblem:
 % unique solution of the principal eigenvalue problem:
\begin{equation}\label{eq:eigenvalue_pb_main}
  \left\{
  \begin{aligned}
    \mcl^{\lambda}_{m,L}\,\varphi&=k^\lambda_{m,L}\,\varphi&\text{in $\Er\times\Theta$},\\
    \nu\cdot\nabla_{\theta}\varphi&=0&\text{over $\Er\times\dr\Theta$},\\
    \varphi&>0&\text{in $\Er\times\Theta$},\\
    \varphi&&\text{is $L$-periodic in $x$}.
  \end{aligned}
  \right.
\end{equation}
%such that
%\[\int_{[0,L]\times\Theta}\varphi_{m,L}^{\lambda}=1.\]
Due to the periodicity of $x$ of the model and to the smoothness of $\dr\Theta$, we are in a compact setting.
Therefore, it follows from the standard Krein-Rutman theory~\cite{KreRut48} that $k^{\lambda}_{m,L}$ is well-defined and is associated to a unique (up to multiplication) principal eigenfunction. %the existence and uniqueness of such a couple $(\varphi^{\lambda}_{m,L},k^{\lambda}_{m,L})\in \mcc^{2}(\Er\times\overline{\Theta})\times\Er$ follows from the standard Krein-Rutman theory~\cite{KreRut48}.
%Several principal eigenvalues/eigenfunctions will be introduced later, and their existence and uniqueness all follow from the Krein-Rutman theory.

%\paragraph{Persistence.}

%\paragraph{Propagation.}

Our first result means that the solutions of the equations of the type \eqref{eq:main} have a spreading speed, and that the spreading speed is given by a formula analogous to the Freidlin-Gärtner formula (see Equation \eqref{eq:freidlin_gartner_intro}). First, we extend the definition of the spreading speed to our framework.
\begin{dfi}\label{dfi:spreading_speed}
  Let~$c>0$. We say that a solution~$u$ of \eqref{eq:main} spreads at speed~$c$ if for all~$c'>c$,
\[\lim_{t\to+\infty}\rho(t,c't)= 0,\]
and for all~$c''\in(0,c)$,
\[\liminf_{t\to+\infty}\rho(t,c''t)>0.\]
We call~$c$ the {spreading speed} of~$u$.
\end{dfi}

\begin{req}
  Using a parabolic Harnack equality (see below), we may infer that~$u$ spreads at speed~$c$ if and only if: for all~$c'>c$,
\[\lim_{t\to+\infty}\pth{\sup_{\theta\in\Theta}u(t,c't,\theta)}= 0,\]
and for all~$c''\in(0,c)$,
\[\liminf_{t\to+\infty}\pth{\inf_{\theta\in\Theta}u(t,c''t,\theta)}>0.\]
\end{req}

As we are interested in the spreading of a population, the least that we should require is that the population does not get extinct. Theorem~2.4 in~\cite{br24} implies that the population persists as long as the mutation coefficient~$m>0$ and the period~$L>0$ satisfy
\[k^0_{m,L}>0.\]
We shall only work with such mutation coefficients and periods.
For example, this condition is satisfied for all $m>0$ and $L>0$ if $r\goq 0$ and $r\not\equiv0$.

We now have all the tools at hand to state the Freidlin-Gärtner formula in our context. 

\begin{thm}[Freidlin-Gärtner formula]\label{thm:freidlin_gartner}
  Take~$m>0$ and~$L>0$ and assume that $u_0\in\mcc^{0}(\Er\times\overline{\Theta})$ is nonnegative, not identically equal to~$0$, bounded and has a support bounded from above.
  \begin{enumerate}[label=$(\roman*)$]
  \item If ${k^0_{m,L}<0}$, then the solution of \eqref{eq:main} gets extinct: $\dabs{\rho(t,\cdot)}_{\infty}\to 0$ as $t\to+\infty$;
  \item If ${k^0_{m,L}>0}$, then the solution of \eqref{eq:main} has a positive spreading speed~$c_{m,L}$ (in the sense of Definition~\ref{dfi:spreading_speed}). The spreading speed~$c_{m,L}$ is given by
   \begin{equation}
    c_{m,L}:=\inf_{\lambda>0}\frac{k^{\lambda}_{m,L}}{\lambda}.
  \end{equation}
  \end{enumerate}
\end{thm}

\begin{req}
  As mentioned in~\cite{br24}, the behaviour of the solution (extinction or persistence) in the borderline case~$k^0_{m,L}=0$ is not known.
\end{req}

\subsection{The spreading speed for different scales}\label{ss:different_scales}

In this section, we shall compute the spreading speed for different scaling limits. We will often go back to the study of equations without the phenotypic dimension. It will thus be convenient to use the notation $c^1(A,R)$ for the spreading speed of a population not structured in phenotype, living in an environment with periodic diffusion coefficient $A\in\mcc^{1,\alpha}(\Er)$ and with periodic intrinsic growth rate $R\in\mcc^{0,\alpha}(\Er)$, with the same period as $A$. 

Let us precisely define $c^1(A,R)$ for two functions $A,R$ which are periodic with a common period $L>0$.
We let $\kappa^{1}(\lambda;A,R)$ be the principal eigenvalue corresponding to the operator~$\mcl^{1,\lambda}$  defined by
\[\mcl^{1,\lambda}\phi:=(A\phi')'-2\lambda A\phi'+(R+\lambda^2A-\lambda A')\phi,\esp\phi\in\mcc^{2,\alpha}(\Er),\]
acting on $L$-periodic functions.
This means that $\kappa^1(\lambda;A,R)$ is the only real number such that there exists a $L$-periodic function $\varphi^{1,\lambda}\in\mcc^{2,\alpha}(\Er)$ satisfying
\[\mcl^{1,\lambda}\varphi^{1,\lambda}=\kappa^1(\lambda;A,R)\varphi^{1,\lambda},\esp\varphi^{1,\lambda}>0.\]
We point out that $\kappa^1(\lambda;A,R)$ does not depend on the common period of $A$ and $R$ that we choose.
Then, we set
\[c^1(A,R):=\inf_{\lambda>0}\frac{\kappa^{1}(\lambda;A,R)}{\lambda},\]
which, by the classical Freidlin-Gärtner formula \eqref{eq:freidlin_gartner_intro}, does indeed correspond to the spreading speed of a population not structured in phenotype, with diffusion coefficient~$A$ and with intrinsic growth rate $R$. 

Unless otherwise specified, we allow $a$ and $\mu$ to depend on both variables $x$ and $\theta$.

\paragraph{As $L\to+\infty$: Possible to forget the phenotype.}

First, we compute the limit of the spreading speed as the period goes to infinity. We assume that $a=a(x)$ is independent of $\theta$ and that $\mu=\mu(\theta)$ is independent of $x$.
For all $x\in\Er$, we let $H_m(x)$ be the principal eigenvalue associated to the \emph{local eigenproblem}: 
\begin{equation}\label{eq:local_problem}
  \left\{
  \begin{aligned}
    m\nabla_{\theta}\cdot\pth{\mu\nabla_{\theta}\psi}+r(x,\cdot)\psi&=H_m(x)\psi&\text{in $\Theta$},\\
    \nu\cdot\nabla_{\theta}\psi&=0&\text{over $\dr\Theta$},\\
    \psi&>0&\text{in $\Theta$}.
  \end{aligned}
  \right.
\end{equation}
%Thanks to Lemma 1 of \cite{BouMir15}, there exists a unique solution $(\psi_x,H(x))\in\mcc^2(\Theta)\times\Er$ to \eqref{}, up to multiplication of $\psi_x$.
We call $H_m(x)$ the \emph{local principal eigenvalue}. As mentioned in the introduction, if $a$ and $\mu$ are constant and $r$ only depends on $\theta$, then $H_m$ is independent of $x$ and the spreading speed of the solution is $2\sqrt{aH_m}$. This is reminiscent of the spreading speed for the Fisher-KPP equation: hence $H_m$ plays the role of a local growth rate.

The following theorem is a consequence of Theorem~\ref{thm:freidlin_gartner} and of the results of~\cite{b23}.

\begin{thm}\label{thm:L_to_inf}
   Fix $m>0$ and assume that $k^{0}_{m,L_0}>0$ for some $L_0>0$. Assume that $a=a(x)$ is independent of $\theta$ and that $\mu=\mu(\theta)$ is independent of $x$. Then $c_{m,L}$ converges as $L\to +\infty$ and
   \[\lim_{L\to +\infty}c_{m,L}=\lim_{L\to+\infty}c^1\pth{a_L,(H_m)_L},\]
   where, with the usual notation, $(H_m)_L(x)=H_m(x/L)$ is the $L$-periodic version of $H_m$.
\end{thm}

Theorem~\ref{thm:L_to_inf} allows one to find a more or less explicit formula for the limit, given in Proposition~\ref{ppn:L_to_inf} below (see also~\cite{HNR11}).
  This formula, in fact, is counter-intuitive. Let us try --~and fail~-- to guess what it should be. For $x\in\Er$, let $c(x)=2\sqrt{a(x)H_m(x)}$ be the spreading speed corresponding to the equation with coefficients \enquote{frozen} at $x$. Since the period is very large, the environment is locally almost homogeneous. Hence, we expect that at position $x$, the instantaneous speed of the invasion is $c(x)$. Averaging over the period, the resulting spreading speed should be the harmonic mean of $x\mapsto c(x)$:
\[c_H:=\pth{\int_0^1\frac{1}{c(y)}\de y}^{-1}.\]
Unfortunately, this is not the case in general. We have indeed, if $a$ is constant and $H$ is not constant:
\begin{equation}\label{eq:cH_vs_cinf}
  c_H<\lim_{L\to+\infty}c_{m,L}.
\end{equation}
This is proved after the proof of Theorem~\ref{thm:L_to_inf}.
The same phenomenon also occurs for the equation without the phenotypic dimension (see~\cite{HNR11}). For bistable equations, the above heuristic does work, see~\cite{DHL22}. 

This counter-intuitive phenomenon has been dubbed the \enquote{tail problem} and is due to the fact that the front is pulled by an extremely small sub-population~\cite{Fre85,Jab12}. This is a true caveat in the modelling properties of equations of the Fisher-KPP type.
A first possibility to overcome it is to  consider the bistable case by adding an Allee effect~\cite{HFR10}.
Biologically, this means that cooperation is needed for propagation, so that a small population cannot survive by itself. 
A second possibility is to first work in finite populations, consider the limit $L\to+\infty$, and only then consider the infinite population limit; with this respect, see~\cite{MRT21}.

\paragraph{As $L\to 0$: Impossible to forget the phenotype.} 

Second, we compute the limit of $c_{m,L}$ as $L\to 0$. Here, we shall not be able to express the limit in terms of $1$-dimensional speeds $c^1(A,R)$.

We introduce a \enquote{homogenised} principal eigenvalue problem. For a function $f:\Er\times\Theta\to\Er$ which is $1$-periodic in $x$, we let
\[\widetilde{f}^A(\theta):=\int_{0}^1f(x,\theta)\de x\]
be the arithmetic mean on a period in the direction $\Er$.  If $f>0$, we let 
\[\widetilde{f}^H(\theta):=\pth{\int_{0}^1\frac{\de x}{f(x,\theta)}}^{-1}\]
be the harmonic mean on a period in the direction $\Er$.
We let $\widetilde{k}_m^{\lambda}$ be the principal eigenvalue associated to the \emph{homogenised eigenproblem}: 
\begin{equation}\label{eq:homo_eig_problem}
  \left\{
  \begin{aligned}
    m\nabla_{\theta}\cdot\pth{\widetilde{\mu}^A\nabla_{\theta}\psi}+\pth{\widetilde{r}^A+\lambda^2\widetilde{a}^H}\psi&=\widetilde{k}^{\lambda}_m\psi&\text{in $\Theta$},\\
    \nu\cdot\nabla_{\theta}\psi&=0&\text{over $\dr\Theta$},\\
    \psi&>0&\text{in $\Theta$}.
  \end{aligned}
  \right.
\end{equation}
We call $\widetilde{k}^{\lambda}_m$ the \emph{homogenised principal eigenvalue}.
The following theorem is a consequence of homogenisation techniques.

\begin{thm}\label{thm:L_to_0}
  Fix $m>0$ and assume that $\widetilde{k}^0_m>0$. 
  Then the limit of $c_{m,L}$ as $L\to 0$ exists and
  \[\lim_{L\to 0}c_{m,L}=\inf_{\lambda>0}\frac{\widetilde{k}^{\lambda}_m}{\lambda}.\]
  Moreover, we have the following variational formula:
  \[\lim_{L\to 0}c_{m,L}=2\,\max_{\psi\in H^1(\Theta),\,\|\psi\|_{L^2(\Theta)}=1,\,R_{\psi}>0}\ \pth{\sqrt{A_{\psi}R_{\psi}}},\]
  where $A_{\psi}:=\displaystyle\int_{\Theta}{\widetilde{a}^H\psi^2}$ and $R_{\psi}:=\displaystyle\int_{\Theta}\pth{\widetilde{r}^A\psi^2-m\widetilde{\mu}^A\abs{\nabla\psi}^2}$.
\end{thm}

  If the coefficients $a$, $\mu$ and $r$ are independent of the phenotype $\theta$, then so are the principal eigenfunctions associated to the homogenised eigenproblem; thus
  \[\widetilde{k}_m^{\lambda}=\widetilde{r}^A+\lambda^2\widetilde{a}^H\]
  and we obtain, by the first equality in Theorem~\ref{thm:L_to_0},
  \[\lim_{L\to 0}c_{m,L}=2\sqrt{\widetilde{a}^H\widetilde{r}^A}.\]
  This is the result proved in~\cite{EHR09} for the equation without the phenotypic variable.

\paragraph{As $m\to 0$ or $m\to+\infty$: Possible to forget the phenotype.} 
Last, we compute the limits of the spreading speed as the mutation scaling parameter $m$ goes to zero or infinity. 
For a function $f$ defined on $\Er\times\Theta$, we let
\[\overline{f}(x):=\frac{1}{\abs{\Theta}}\int_{\Theta}f(x,\theta)\de\theta\]
be the arithmetic mean of $f$ in the direction $\Theta$. Note that $\overline{f}$ differs from the notation $\widetilde{f}^A$, which is the arithmetic mean of $f$ on a period in the direction $\Er$.

\begin{thm}\label{thm:mutation}
  Fix $L>0$. 
  \begin{enumerate}[label=$(\roman*)$]
  \item %Assume that there exist $\delta>0$ and $\underline{m}>0$ such that $k^0_{m,L}>\delta$ for all $m>\underline{m}$.
    Assume that $a=a(x)$ and $\mu=\mu(x)$ are independent of $\theta$. 
    Assume moreover that $\kappa^1(0;a,\overline{r})>0$.
    Then the limit of $c_{m,L}$ as $m\to+\infty$ exists and
    \begin{equation*}
      \lim_{m\to +\infty}c_{m,L}=c^1\pth{a,\overline{r}}.
    \end{equation*}
  \item %Assume that there exist $\delta>0$ and $\overline{m}>0$ such that $k^0_{m,L}>\delta$ for all $m\in\ioo{0,\overline{m}}$. 
    Assume that $a=a(x)$ is independent of $\theta$ and that $\mu=\mu(\theta)$ is independent of~$x$.
    Assume moreover that there exists $\theta\in\overline{\Theta}$ such that $\kappa^1(0;\, a(\cdot),\, r(\cdot,\theta))>0$.
    Then the limit of $c_{m,L}$ as $m\to 0$ exists and
     \begin{equation*}
      \lim_{m\to 0}c_{m,L}=\max_{\theta\in\overline{\Theta}\ / \ \kappa^1(0;\, a(\cdot),\, r(\cdot,\theta))>0}c^1(a(\cdot),r(\cdot,\theta)).
     \end{equation*}
%     Here, if $\kappa^1(0;a(\cdot),r(\cdot,\sigma))<0$ for some $\sigma\in\overline{\Theta}$, we understand $c^1(a(\cdot),r(\cdot,\sigma)):=-\infty$.
   \item If the conditions for the first point and second points are satisfied, then 
     \[\lim_{m\to +\infty}c_{m,L}\loq\lim_{m\to 0}c_{m,L}.\]
  \end{enumerate}
\end{thm}

In this theorem, the assumptions that $\kappa^1(0;\, a,\,\overline{r})>0$ and that $\kappa^1(0;\,a(\cdot),\, r(\cdot,\theta))>0$ for some $\theta\in\overline{\Theta}$ ensure that in the limiting regime, there is indeed a spreading speed.
Item $(i)$ of the theorem is a consequence of homogenisation techniques, and item $(ii)$ is a consequence of the results in~\cite{b23}.

\subsection{The effect of the mutations}\label{ss:mutation_effect}

Using the previous technical mathematical results, we are ready to give several consequences which have an easy biological interpretation and which highlight the role of the phenotypic dimension in the qualitative behaviour of the population.
We assume throughout Subsection~\ref{ss:mutation_effect} that $a=a(x)$ is independent of $\theta$ and that $\mu=\mu(\theta)$ is independent of $x$. 
We define 
\[\gamma:=\lim_{m\to 0}\lim_{L\to+\infty}c_{m,L}-\lim_{L\to +\infty}\lim_{m\to 0}c_{m,L}.\]
We assume that there is no phenotype which is the best everywhere:
  \begin{equation}\label{eq:assumption_hetero}
    \forall\theta\in\overline{\Theta},\ \ \exists(x,\,\sigma)\in\Er\times\overline{\Theta},\esp r(x,\sigma)>r(x,\theta).
  \end{equation}
  The following theorem shows that 
  that under assumption~\eqref{eq:assumption_hetero},
  $\gamma$ is well-defined,
  is positive and can be interpreted as the maximal potential effect of mutations on the spreading speed.

\begin{thm}\label{thm:bio_csq} 
  We assume that $a$ is independent of $\theta$ and that $\mu$ is independent of $x$.
  We also assume that there exist $\underline{L}>0$ and $\theta\in\overline{\Theta}$ such that $\kappa^1(0;a,r_{\underline{L}}(\cdot,\theta))>0$.
  Finally, we assume that \eqref{eq:assumption_hetero} holds.
  Then, there exists $L_0>0$ such that for all $L>L_0$, there exists an optimal mutation coefficient $m^*(L)>0$, in the sense that
  \[c_{m^*(L),L}=\sup_{m>0}c_{m,L}.\]
  Moreover, $\gamma$ is well-defined, is positive, and for all $\eps>0$, there exists $L^{\eps}\goq L_0$ such that for all $L>L^{\eps}$,
  \[c_{m^*(L),L}>\lim_{m\to 0}c_{m,L}+\gamma-\eps\] %\goq \lim_{m\to +\infty}c_{m,L}+\gamma-\eps.\]
\end{thm}

Let us comment on this theorem.

\medskip

First, the fact that there is an optimal mutation coefficient $m^*(L)>0$, together with the inequalities, suggest that there is a trade-off between two phenomena: too small a mutation coefficient implies that the population cannot benefit from the diversity of phenotypes; too high a mutation coefficient implies that the population suffers from a mutation load, \ie suffers from the influence of the least adapted phenotypes.
  As the following proposition shows, this trade-off disappears in the limits $L\to+\infty$ and $L\to 0$: in both cases, the environment is \enquote{almost homogeneous}, so that having a higher capacity of adaptation is useless, but increases the mutation load.
\begin{ppn}\label{ppn:limit_dec_mu}
  We assume that $a$ is independent of $\theta$ and that $\mu$ is independent of $x$. We let $I\subset\ioo{0,+\infty}$ be an interval.
  \begin{enumerate}[label=$(\roman*)$]
  \item Assume that $\widetilde{k}^0_m>0$ for all $m\in I$. %, for all $L\in\ioo{0,\overline{L}}$, we have $k^0_{m,L}>\delta$.
    Then the function
    \[m\in I\mapsto\lim_{L\to 0}c_{m,L}\]
    is well-defined and nonincreasing.
  \item Assume that there exists $L_0>0$ such that for all $m\in I$, we have $k^0_{m,L_0}>0$.
    Then the function
    \[m\in I\mapsto\lim_{L\to +\infty}c_{m,L}\]
  is well-defined and nonincreasing.
  \end{enumerate}
  Moreover, if there exists $x\in\Er$ such that $\theta\mapsto r(x,\theta)$ is not constant, then we can replace \enquote{nonincreasing} by \enquote{strictly decreasing} in both statements.
\end{ppn}

Second, the fact that $\gamma>0$ means that for an \enquote{infinite period}, very small mutations have a non-vanishing effect on the spreading speed. The value $\gamma$ describes the maximal effect on the speed which mutations can have. With this respect, see also~\cite{EllCor12,MBC19,rbzp24}. 

\medskip

Third, we state a slightly counter-intuitive result.
%increasing the ability of persistence is not equivalent to increasing the spreading speed.

\begin{ppn}\label{ppn:persistence_decr}
  Fix $L>0$. The function $m\mapsto k^{0}_{m,L}$ is nonincreasing.
\end{ppn}

Theorem~\ref{thm:bio_csq} implies that for $L$ large enough, $m\mapsto c_{m,L}$ is not nonincreasing. But, as mentioned above, the ability of persistence can be measured by the value of $k^0_{m,L}$. Thus, Proposition~\ref{ppn:persistence_decr} means that \emph{making the persistence easier is not equivalent to increasing the spreading speed}.

\medskip

Last, we state a proposition which brings out the different qualitative properties of the two phenomena, spreading and mutation. For the period, there is no trade-off: the higher the period, the higher the speed. This property also holds when the phenotype is not taken into account~\cite{Nad10}. 

\begin{ppn}\label{ppn:inter_L} 
  Fix $m>0$. Let $L_0>0$ such that $k^{0}_{m,L_0}>0$.
  The function $L\mapsto c_{m,L}$ is positive, continuous and nondecreasing on $\ifo{L_0,+\infty}$.
\end{ppn}

In Section~\ref{s:fg}, we will show that solutions of equations of the type \eqref{eq:main} have a spreading speed which is given by a formula analogous to the Freidlin-Gärtner formula (Theorem~\ref{thm:freidlin_gartner}). Using this formula, we will then look at the behaviour of the spreading speed of the solution of \eqref{eq:main} as $L\to 0$, $L\to +\infty$ (Section~\ref{s:period}), $m\to 0$ and $m\to +\infty$ (Section~\ref{s:mutation}). The qualitative properties about the effect of the mutations are proved in Subsections~\ref{ss:inter_period},~\ref{ss:proof_bio} and~\ref{ss:other}.

\section{The Freidlin-Gärtner formula (Theorem~\ref{thm:freidlin_gartner})}\label{s:fg}
In this part, we prove Theorem~\ref{thm:freidlin_gartner}. For notational convenience, we shall assume throughout that $L=1$ and $m=1$. 
The notations thus simplify:
\[c=c_{m,L},\esp k^{\lambda}=k^{\lambda}_{m,L}.\]

\subsection{A Harnack inequality}

We first state a nonstandard version of the parabolic Harnack inequality,
which was found by Alfaro, Berestycki and Raoul in~\cite{ABR17} and which is well adapted to our context.
Contrarily to the standard parabolic Harnack inequality, in this new inequality, there is no delay in time.

\begin{thm}[\cite{ABR17}]\label{thm:harnack_abr}
  %Let $K>0$ and $f\in L^{\infty}(\ifo{0,+\infty}\times \Er\times\Theta)$ such that $\dabs{f}_{\infty}<K$, $\dabs{a}_{\infty}<K$ and $\dabs{\dr_xa}_{\infty}<K$.
  Let  $f\in L^{\infty}(\ifo{0,+\infty}\times \Er\times\Theta)$.
  Let $\delta>0$ and let $u(t,x,\theta)$ be a bounded nonnegative weak solution of
\begin{equation}\label{eq:harnack_abr_pheno}
  \left\{
  \begin{aligned}
    \dr_tu(t,x,\theta)&=\dr_x(a\dr_xu)+\nabla_{\theta}\cdot(\mu\nabla_{\theta}u)+f(t,x,\theta)u,&t>0\vg (x,\theta)\in\Er\times\Theta,\\
    \nu\cdot\nabla_{\theta} u(t,x,\theta)&=0,&(x,\theta)\in\Er\times\dr\Theta.
  \end{aligned}
  \right.
\end{equation}
%Then there exists a constant $C>0$ depending only on $\dabs{u}_{\infty}$, $\delta$, $K$, $a$, $\mu$ and $\Theta$ such that:
Then there exists a constant $C>0$ depending only on $\delta$, $\dabs{u}_{\infty}$, $\dabs{f}_{\infty}$, $a$, $\mu$ and~$\Theta$ such that:
  \[\forall t_0\goq 1,\, \forall x_0\in\Er,\qquad\sup_{\theta\in\Theta}u(t_0,x_0,\theta)\loq C\inf_{\theta\in\Theta}u(t_0,x_0,\theta)+\delta.\]

\end{thm}

%\begin{thm}[\cite{ABR17}]\label{thm:harnack_abr}
%  Let $K>0$ and $f\in L^{\infty}(\ifo{0,+\infty}\times \Er\times\Theta)$ such that $\dabs{f}_{\infty}<K$, $\dabs{a}_{\infty}<K$ and $\dabs{\dr_xa}_{\infty}<K$.
%  Set $U>0$ and $\delta>0$.
%
%   Then there exists a constant $C>0$ depending only on $K$, $U$, $\delta$, $a$, $\mu$ and $\Theta$ such that if the nonnegative weak solutions $u(t,x,\theta)$ of
%%  Then there exists a constant $C>0$ depending only on $K$, $U$, $\delta$, $\eta$, $P$ and $diam(\Theta)$ such that a nonnegative function $u(t,x,\theta)$ which satisfies $\dabs{u}_{\infty}\loq U$ and
%\begin{equation}\label{eq:harnack_abr_pheno}
%  \left\{
%  \begin{aligned}
%    \dr_tu(t,x,\theta)&=\dr_x(a\dr_xu)+\nabla_{\theta}\cdot(\mu\nabla_{\theta}u)+f(t,x,\theta)u,&t>0\vg (x,\theta)\in\Er\times\Theta,\\
%    u(0,x,\theta)&=u_0(x,\theta),&(x,\theta)\in\Er\times\Theta,\\
%    \nu\cdot\nabla_{\theta} u(t,x,\theta)&=0,&(x,\theta)\in\Er\times\dr\Theta,
%  \end{aligned}
%  \right.
%\end{equation}
%satisfies $\dabs{u}_{\infty}\loq U$, then $u$ 
%satisfies also:
%  \[\forall t_0\goq 1,\, \forall x_0\in\Er,\qquad\sup_{\theta\in\Theta}u(t_0,x_0,\theta)\loq C\inf_{\theta\in\Theta}u(t_0,x_0,\theta)+\delta.\]
%
%\end{thm}

\begin{proof}
  The proof is almost the same as that of~\cite[Theorem 2.6]{ABR17}, which holds only in the full space. The only difference is that we need to deal with the boundary $\dr\Theta$.
  
  Without loss of generality, we prove the result for $t_0=2$ and $x_0=0$. We consider $R>diam(\Theta)$.
  For convenience, we note
  \[U:=\dabs{u}_{\infty},\qquad K:=\max\pth{\dabs{f}_{\infty},\, \dabs{a}_{\infty},\,\dabs{\dr_xa}_{\infty}}.\]
  As in~\cite{ABR17}, we introduce the following function, defined on $\ifo{0,+\infty}\times\Er\times\Theta$ and independent of the variable $\theta$:
  \[\phi(t,x,\theta):=e^{K(t-1)}\cro{\max_{\abs{x'}\loq\beta R\vg \theta'\in\Theta}u(1,x',\theta')+\frac{2U}{(\beta R)^2}K(1+\beta R)(t-1)+\frac{Ux^2}{(\beta R)^2}},\]
  for some $\beta>1$ to be determined later. We have:
  \begin{align*}
    &\dr_t\phi(t,x,\theta)-\dr_x(a\dr_x\phi)-\nabla_{\theta}\cdot(\mu\nabla_{\theta}\phi)-f(t,x,\theta)\phi\\
    &\esp \goq \pth{K\phi+e^{K(t-1)}\times\frac{2U}{(\beta R)^2}K(1+\beta R)}-e^{K(t-1)}\times \frac{2U}{(\beta R)^2}\pth{a+x\dr_xa}-K\phi.
  \end{align*}
  Therefore, since $\dabs{a}_{\infty}\loq K$ and $\dabs{\dr_xa}_{\infty}\loq K$, we have
  \[\dr_t\phi(t,x,\theta)-\dr_x(a\dr_x\phi)-\nabla_{\theta}\cdot(\mu\nabla_{\theta}\phi)-f(t,x,\theta)\phi\goq 0\]
  on $\ioo{1,2}\times B(0,\beta R)\times\Theta$.
  We also have the boundary conditions: $\phi(1,x,\theta)\goq u(1,x,\theta)$ on $B(0,\beta R)\times\Theta$ and $\phi(t,x,\theta)\goq u(t,x,\theta)$ on $(\dr B(0,\beta R))\times\Theta$.
  Finally, we have $\nu\cdot\nabla_{\theta}\phi=0$ and $\nu\cdot \nabla_{\theta} u=0$ on $B(0,\beta R)\times \dr\Theta$.
  Therefore we can apply the comparison principle:
  \[\phi\goq u\comment{on $\cro{1,2}\times B(0,\beta R)\times\Theta$}.\]
  In particular, for $\beta>0$ large enough (depending only on $R$, $U$, $K$ and $\delta$),
  \begin{align*}
    \sup_{\abs{x}<R\vg \theta\in\Theta}u(2,x,\theta)&\loq \sup_{\abs{x}<R\vg \theta\in\Theta}\phi(2,x,\theta)\\
    &\loq e^K\sup_{\abs{x}<\beta R\vg \theta\in\Theta}u(1,x,\theta)+\delta.
  \end{align*}
  We now wish to apply the parabolic Harnack inequality up to the boundary of $\Theta$.
  In order to do this, we use the regularity of $\dr\Theta$ and the Neumann boundary condition to extend the solution $u$ to the solution of a uniformly parabolic equation on a larger domain $\Er\times\widetilde{\Theta}$ with $\overline{\Theta}\subset\widetilde{\Theta}$.
  See \eg the argument for the proof of Step~2 in~\cite[Theorem 3.1]{BerRos09}.
  Applying the classical interior parabolic Harnack inequality in $\Er\times\widetilde{\Theta}$
  implies, in our context,
  that we may apply the parabolic Harnack inequality up to the boundary of $\Theta$: there exists a constant $C>0$, depending only on $K$, $\beta R$, $a$, $\mu$ and $\Theta$, such that
  \[e^K\sup_{\abs{x}<\beta R\vg \theta\in\Theta}u(1,x,\theta)\loq C\inf_{\abs{x}<\beta R\vg \theta\in\Theta}u(2,x,\theta).\]
  This gives:
  \[\sup_{\abs{x}<R\vg \theta\in\Theta}u(2,x,\theta)\loq C\inf_{\abs{x}<\beta R\vg \theta\in\Theta}u(2,x,\theta)+\delta.\]
  This implies
  \[\sup_{\theta\in\Theta}u(2,0,\theta)\loq C\inf_{\theta\in\Theta}u(2,0,\theta)+\delta.\]
  This is what we wanted to prove

\end{proof}

\subsection{An estimate on the principal eigenvalue}

\begin{lem}\label{lem:estimate_klambda}
  For all $\lambda>0$,
  \[\lambda^2\inf a-\lambda\dabs{\dr_xa}_{\infty}-\dabs{r}_{\infty}\loq k^{\lambda}\loq \lambda^2\sup a+\lambda\dabs{\dr_xa}_{\infty}+\dabs{r}_{\infty}.\]
\end{lem}

\begin{proof}
  For a $1$-periodic in $x$ function $g\in\mcc^{0,\alpha}(\Er\times\overline{\Theta})$, we let $\kappa(g)$ be the principal eigenvalue associated to the operator defined by
  \[\mcl_g u=\dr_x(a\dr_xu)+\nabla_{\theta}\cdot\pth{\mu\nabla_{\theta} u}-2\lambda a\dr_xu+gu,\]
  with periodicity in $x$ and with Neumann boundary conditions on $\Er\times\dr\Theta$. By definition, we have
  \[k^{\lambda}=\kappa(\lambda^2a-\lambda\dr_xa+r).\]
  Moreover, if $f_1\loq f_2$ then $\kappa(f_1)\loq \kappa(f_2)$. Therefore,
  \[\kappa(\lambda^2\inf a-\lambda\dabs{\dr_xa}_{\infty}-\dabs{r}_{\infty})\loq k^{\lambda}\loq \kappa(\lambda^2\sup a+\lambda\dabs{\dr_xa}_{\infty}+\dabs{r}_{\infty}).\]
  Finally, $\varphi\equiv 1$ is an eigenfunction of the operator $\mcl_g$ when $g$ is constant, which implies:
  \begin{align*}
    k(\lambda^2\inf a-\lambda\dabs{\dr_xa}_{\infty}-\dabs{r}_{\infty})&=\lambda^2\inf a-\lambda\dabs{\dr_xa}_{\infty}-\dabs{r}_{\infty},\\
    k(\lambda^2\sup a+\lambda\dabs{\dr_xa}_{\infty}+\dabs{r}_{\infty})&=\lambda^2\sup a+\lambda\dabs{\dr_xa}_{\infty}+\dabs{r}_{\infty}.
  \end{align*}
  The conclusion follows.
\end{proof}

\subsection{Proof of Theorem~\ref{thm:freidlin_gartner}}

The extinction case for $k^0<0$ is proved in~\cite{br24}.
From now on, we assume that $k^0>0$.
We define
\[c_{FG}:=\inf_{\lambda>0}\frac{k^{\lambda}}{\lambda}.\]
First, $\mcl^{\lambda}$ is the adjoint of $\mcl^{-\lambda}$, so the function $\lambda\in\Er\mapsto k^{\lambda}$ is symmetric.
Further, the argument of the proof of~\cite[Lemma~3.1]{BHR05-2} implies that $\lambda\in\Er\mapsto k^{\lambda}$ is also convex.
Hence $k^{\lambda}>0$ for all $\lambda\in\Er$.
Therefore, the lower bound in~Lemma~\ref{lem:estimate_klambda} implies that $c_{FG}>0$.

   Now, let us show that a solution $u$ of \eqref{eq:main} starting from a bounded nonnegative initial condition with support bounded from above,
   has a spreading speed $c_{FG}$.
   The proof is cut into two parts:
   Proposition~\ref{ppn:logcfg} below means that the first requirement of Definition~\ref{dfi:spreading_speed} is satisfied;
   Proposition~\ref{ppn:goqcfg} below means that the second requirement of Definition~\ref{dfi:spreading_speed} is satisfied.

\begin{ppn}\label{ppn:logcfg}
  Assume that $u_0\in\mcc^{0}(\Er\times\overline{\Theta})$ is nonnegative, bounded and has a support bounded from above. Let $u$ be the solution of \eqref{eq:main}. For each $c'>c_{FG}$,
  \[\lim_{t\to+\infty}\rho(t,c't)=0.\]
\end{ppn}

\begin{proof}
  Let $c'>c_{FG}$ and let $\lambda>0$ such that $c'':=k^{\lambda}/\lambda$ satisfies $c_{FG}\loq c''<c'$. Let $\varphi^{\lambda}>0$ be a (periodic) principal eigenfunction of $\mcl^{\lambda}$.
  Let \[\bar{u}(t,x,\theta):=Ce^{-\lambda(x-c''t)}\varphi^{\lambda}(x,\theta),\]
  for some constant $C>0$ to be chosen later.
  Short computations show that $\bar{u}$ solves \[\dr_t\bar{u}=\dr_x(a\dr_x \bar{u})+\nabla_{\theta}\cdot\pth{\mu\nabla_{\theta}\bar{u}}+r(x,\theta)\bar{u}.\]
  Since $u_0$ is bounded and has a support bounded from above, there exists $C>0$ so large that $u_0\loq \bar{u}(0,\cdot,\cdot)$. We obtain by the parabolic comparison principle: for all $t\goq 0$,
  \begin{equation}\label{eq:comp_u_ubar}
  u(t,\cdot,\cdot)\loq \bar{u}(t,\cdot,\cdot).    
  \end{equation}
  Since $\varphi^{\lambda}$ is continuous and periodic, $\varphi^{\lambda}$ is bounded. Thus for $c'>c''$, $\bar{u}(t,c't,\theta)$ converges to~$0$ uniformly in $(x,\,\theta)$ as $t\to+\infty$. With \eqref{eq:comp_u_ubar}, this proves the proposition.
\end{proof}

\begin{ppn}\label{ppn:goqcfg}
  Assume that $u_0\in\mcc^{0}(\Er\times\overline{\Theta})$ is nonnegative, bounded and has a support bounded from above. Let $u$ be the solution of \eqref{eq:main}. For each $c''\in(0,c_{FG})$,
  \[\liminf_{t\to+\infty}\rho(t,c''t)>0.\]
\end{ppn}

In order to prove the proposition, we first consider the auxiliary intrinsic growth rates defined for $\eps>0$ by
\[r^\eps(x,\theta):=r(x,\theta)-\eps,\]
and the following auxiliary \emph{local} problems:
\begin{equation}\label{eq:aux_uepsilon}
  \left\{
  \begin{aligned}
    \dr_tu^{\eps}(t,x,\theta)&=\dr_x(a\dr_xu^{\eps})+\nabla_{\theta}\cdot\pth{\mu\nabla_{\theta}u^{\eps}}+u^{\eps}\pth{r^{\eps}(x,\theta)-G u^{\eps}},&t> 1\vg(x,\theta)\in\Er\times\Theta,\\
    u^{\eps}(1,x,\theta)&=\frac{1}{2}u(1,x,\theta),&(x,\theta)\in\Er\times\Theta,\\
    \nu\cdot\nabla_{\theta} u^{\eps}(t,x,\theta)&=0,&t\goq 1\vg(x,\theta)\in\Er\times\dr\Theta,
  \end{aligned}
  \right.
\end{equation}
where the constant $G>0$ is to be chosen in the proof of Lemma~\ref{lem:estimation_u_below}.
We note that $0<u^{\eps}(1,\cdot,\cdot)<u^{\eps}(1,\cdot,\cdot)$ thanks to the strong parabolic maximum principle. We will first show that for $G>0$ well chosen, $u$ spreads at least as fast as $u^{\eps}$. Then, we will estimate the spreading speed of $u^{\eps}$. Importantly, let us point out that $u^{\eps}$ spreads at a speed \emph{independent} of $G$ (the speed, indeed, only depends on the linear problem).

\begin{lem}\label{lem:estimation_u_below}
    Let $u$ be the solution of the nonlocal problem \eqref{eq:main}. Let $\eps>0$. There exists $G>0$ such that the solution $u^\eps$ of \eqref{eq:aux_uepsilon} satisfies, for all $t\goq 1$, $x\in\Er$, $\theta\in\Theta$,
    \[u(t,x,\theta)\goq u^{\eps}(t,x,\theta).\]
\end{lem}

\begin{proof}
  Let us write the first line of the nonlocal problem \eqref{eq:main} as
  \[
  \dr_tu(t,x,\theta)-\dr_x(a\dr_xu)-\nabla_{\theta}\cdot\pth{\mu\nabla_{\theta}u}=f(t,x,\theta)u
  \]
  with
  \[f(t,x,\theta)=r(x,\theta)-\rho(t,x).\]
%  By Lemma~\ref{lem:apriori}, and since $r$ is bounded, there exist $U>0$ and $K'>0$ such that
  %  \[0\loq u\loq U,\esp -K'\loq f(x,\theta)\loq K'.\]
  By the existence theorem of~\cite{br24}, $u$ and $\rho$ are globally bounded.
  Further, $r$ is bounded, so $f$ is also globally bounded.
  Therefore, we can apply Theorem~\ref{thm:harnack_abr}: for each $\delta>0$, there exists $C_{\delta}>0$ such that for all $t\goq 1$, $x\in\Er$, $\theta\in \Theta$,
  \[C_{\delta}u(t,x,\theta)+\delta\goq \sup_{\sigma\in \Theta}u(t,x,\sigma).\]
  This implies that for each $t\goq 1$, $x\in\Er$, $\theta\in \Theta$,
  \[\rho(t,x)=\int_{\Theta}u(t,x,\sigma)\de\sigma\loq \abs{\Theta}C_{\delta}u(t,x,\theta)+\abs{\Theta}\delta,\]
  where $\abs{\Theta}$ denotes the volume of $\Theta$.
  We choose $\delta$ small enough that $\abs{\Theta}\delta<\eps/2$:
  \begin{equation}\label{eq:estimate_rho_above}
    \rho(t,x) \loq \abs{\Theta}C_{\delta}u(t,x,\theta)+\eps/2.
  \end{equation}
  Now, we set $G:=\abs{\Theta}C_{\delta}$ and we show that for this value of $G$, the lemma holds.
  Let
  \[\mct:=\acc{t\goq 1 \,/\, \exists (x_0,\theta_0)\in\Er\times\Theta,\, u(t,x_0,\theta_0)=u^{\eps}(t,x_0,\theta_0)}\]
  If $\mct=\emptyset$, then we are done.
  Now, we assume that we have $\mct\neq\emptyset$ and we let $t_0:=\inf\mct$.
  By continuity of $u$ and $u^{\eps}$,
  we have $t_0\in\mct$. Therefore, there exists $(x_0,\theta_0)\in\Er\times\Theta$
  such that
  \[u(t_0,x_0,\theta_0)=u^{\eps}(t_0,x_0,\theta_0).\]
  We recall that $u^{\eps}(1,\cdot,\cdot)<u(1,\cdot,\cdot)$.
  Then we must have $t_0>1$ and
  \[\cro{\dr_t-\pth{\dr_x(a\dr_x)+\nabla_{\theta}\cdot\pth{\mu\nabla_{\theta}}}}(u-u^{\eps})(t_0,x_0,\theta_0)\loq 0.\]
  This implies:
  \[\pth{r(x_0,\theta_0)-\rho(t_0,x_0)}u(t_0,x_0,\theta_0)-\pth{r^{\eps}(x_0,\theta_0)-G u^{\eps}(t_0,x_0,\theta_0)}u^{\eps}(t_0,x_0,\theta_0)\loq 0.\]
  Recalling that $u(t_0,x_0,\theta_0)=u^{\eps}(t_0,x_0,\theta_0)>0$, we get
  \begin{equation*}
    r(x_0,\theta_0)-\rho(t_0,x_0)\loq r^{\eps}(x_0,\theta_0)-G u(t_0,x_0,\theta_0).
  \end{equation*}
  Finally, we note that $r-r^{\eps}=\eps$. 
  This yields the contact condition:
  \begin{equation}
%    \left\{
    \begin{aligned}
      G u(t_0,x_0,\theta_0)&\loq \rho(t_0,x_0)-\eps.%&\text{if $\theta\in\Theta$},\\
%      G u(t_0,x_0,\theta_0)&\loq \rho(t_0,x_0)-3\eps'-\bar{\rho}&\text{if  $\theta\in\Theta-\Theta$}.
    \end{aligned}
%    \right.
  \end{equation}
  This contact condition is impossible by virtue of \eqref{eq:estimate_rho_above}
  and of the definition of $G$. 
  Therefore, we must have $\mct=\emptyset$. 

\end{proof}

\begin{proof}[Proof of Proposition~\ref{ppn:goqcfg}]
  Let $k^{\eps,\lambda}$ be the principal eigenvalue of the operator
  \[\mcl^{\eps,\lambda}:=\dr_x(a\dr_x)+\nabla_{\theta}\cdot\pth{\mu\nabla_{\theta}}-2\lambda a\dr_x+(r^{\eps}(x,\theta)+\lambda^2a+\dr_xa),\]
  acting on $\mcc^2(\Er\times\overline{\Theta})$ functions which are $1$-periodic in $x$; that is, $k^{\eps,\lambda}$ is the only value such that there exists $\varphi^{\eps,\lambda}\in\mcc^2(\Er\times\overline{\Theta})$ satisfying
  \begin{equation}
  \left\{
  \begin{aligned}
    \mcl^{\eps,\lambda}\varphi^{\eps,\lambda}&={k}^{\eps,\lambda}\varphi^{\eps,\lambda}&\text{in $\Er\times\Theta$},\\
    \nu\cdot\nabla_{\theta}\varphi^{\eps,\lambda}&=0&\text{over $\Er\times\dr\Theta$},\\
    \varphi^{\eps,\lambda}&>0&\text{in $\Er\times\Theta$},\\
    \varphi^{\eps,\lambda}&&\comment{$1$-periodic in $x$}.
  \end{aligned}
  \right.
\end{equation}
  Let us denote by
  \[c^{\eps}_{FG}:=\inf_{\lambda>0}\frac{k^{\eps,\lambda}}{\lambda}\] the Freidlin-Gärtner speed corresponding to $u^{\eps}$.
  By Section 8 of~\cite{Wei02}, the solution $u^{\eps}$ of the local problem \eqref{eq:aux_uepsilon} spreads (in the definition of~\cite{Wei02}) at speed $c^{\eps}_{FG}$. This implies in particular that for all $c''\in(0,c_{FG}^{\eps})$,
  \[\lim_{t\to+\infty}\pth{\inf_{\theta\in\Theta}u^{\eps}(t,c''t,\theta)}>0.\]
  Since $r^{\eps}=r-\eps$, we have $k^{\eps,\lambda}=k^{\lambda}-\eps$.
  Since $c_{FG}>0$, we have $\displaystyle \inf_{\lambda\goq 0}k^{\lambda}>0$.
  Thus:
  \[c^{\eps}_{FG}=\inf_{\lambda>0}\frac{k^{\lambda}-\eps}{\lambda}=c_{FG}+o_{\eps\to 0}(1).\]
  Let $c''\in(0,c_{FG})$.
There exists $\eps>0$ so small that $c''<c^{\eps}_{FG}$. We then have:
\begin{align*}
  \liminf_{t\to+\infty}\rho(t,c''t)&\goq \liminf_{t\to+\infty}\int_{\Theta}u^{\eps}(t,c''t,\sigma)\de\sigma\goq\abs{\Theta}\lim_{t\to+\infty}\pth{\inf_{\theta\in\Theta}u^{\eps}(t,c''t,\theta)}>0.
\end{align*}
This proves the proposition. 

\end{proof}

\section{Behaviour of the spreading speed with respect to the period}\label{s:period}

This part is devoted to the proof of Proposition~\ref{ppn:inter_L} and Theorems~\ref{thm:L_to_inf} and~\ref{thm:L_to_0}.
For notational convenience, we shall assume throughout that $m=1$. 
The notations thus simplify:
\[c_L=c_{m,L},\esp k^{\lambda}_L=k^{\lambda}_{m,L}.\]
In Subsection~\ref{ss:inter_period}, we prove that the function $L\mapsto c_L$ is continuous and nondecreasing (Proposition~\ref{ppn:inter_L}).
In Subsection~\ref{ss:period_inf}, we compute the limit of $c_L$ as $L\to+\infty$ (Theorem~\ref{thm:L_to_inf}).
In Subsection~\ref{ss:period_0}, we compute the limit of $c_L$ as $L\to 0$ (Theorem~\ref{thm:L_to_0}).

\subsection{Continuity and monotony (Proposition~\ref{ppn:inter_L})}\label{ss:inter_period}

\begin{proof}[Proof of Proposition~\ref{ppn:inter_L}, part 1]
  Take $L_0>0$ such that $k^0_{L_0}>\delta$.
  Take $L_1>L_0$. Using Proposition~\ref{ppn:LKL_nondecreasing} below, we have $k^0_{L}>\delta$ for all $L\in\cro{L_0,L_1}$, so $c_L$ is well defined and positive for all $L\in\cro{L_0,L_1}$.
  We prove that the function $L\mapsto c_L$ is continuous on $\cro{L_0,L_1}$.
  The function $(\lambda,L)\mapsto k_L^{\lambda}$ is continuous.
  Hence, as $\lambda\to 0$, we have ${k_L^{\lambda}}/{\lambda}\to+\infty$ uniformly in $L\in\cro{L_0,L_1}$.
  Moreover, by Lemma~\ref{lem:estimate_klambda}, as $\lambda\to+\infty$, we have ${k_L^{\lambda}}/{\lambda}\to+\infty$ uniformly in $L\in\cro{L_0,L_1}$.
  Thus, there exist $\lambda_1>0$ and $\lambda_2>\lambda_1$ such that for all $L\in\cro{L_0,L_1}$,
  \[c_L=\min_{\lambda\in\cro{\lambda_1,\lambda_2}}\frac{k_L^{\lambda}}{\lambda}.\]
  This implies that the function $L\mapsto c_L$ is continuous on $\ioo{L_0,L_1}$, thus on $(L_0,+\infty)$.
\end{proof}

The second part of Proposition~\ref{ppn:inter_L}, \ie the fact that $L\mapsto c_L$ is nondecreasing, is a consequence of the following proposition and of Theorem~\ref{thm:freidlin_gartner}.

\begin{ppn}\label{ppn:LKL_nondecreasing}
  For all $\lambda\in\Er$, the function $L\mapsto k_L^{\lambda}$ is nondecreasing on $\ioo{0,+\infty}$. 
\end{ppn}

The proof of the proposition is almost the same as that of~\cite[Proposition 4.1]{Nad10}. We give an idea of the proof to focus on the assumptions that are needed.
%As in~\cite{Nad10}, we will use several transformations of the parameters of the elliptic equations.
For a $1$-periodic function $f:\Er\times\Theta\to\Er$, we denote by $f_L(x,\theta):=f\pth{\frac{x}{L},\theta}$ its $L$-periodic version.

Take functions $q\in \mcc^{0,\alpha}(\Er\times\overline{\Theta})$ and $g\in \mcc^{0,\alpha}(\Er\times\overline{\Theta})$ which are $1$-periodic in $x$. We denote by $\chi[q,g]$ the principal eigenvalue associated to the eigenvalue problem
\begin{equation*}
  \left\{
  \begin{aligned}
    \dr_x(a_L\dr_x\varphi)&+\nabla_{\theta}\cdot\pth{\mu_L\nabla_{\theta}\varphi}-q_L\dr_x\varphi +g_L\varphi=\chi[q,g]\varphi&\text{in $\Er\times\Theta$},\\
    \nu\cdot\nabla_{\theta}\varphi&=0&\text{over $\Er\times\dr\Theta$},\\
    \varphi&&\text{is $1$-periodic in $x$},\\
    \varphi&>0&\text{in $\Er\times\Theta$}.
  \end{aligned}
  \right.
\end{equation*}
%In particular, $k_L^{\lambda}=\chi_{\lambda}[L,0,r].$
%Take functions $q\in \mcc^{0,\alpha}(\Er\times\overline{\Theta})$ and $g\in \mcc^{0,\alpha}(\Er\times\overline{\Theta})$ which are $1$-periodic in $x$. For $\lambda>0$, we denote by $\chi_{\lambda}[L,q,g]$ the principal eigenvalue associated to the eigenvalue problem
%\begin{equation*}
%  \left\{
%  \begin{aligned}
%    \dr_x(a_L\dr_x\varphi)&+\nabla_{\theta}\cdot\pth{\mu_L\nabla_{\theta}\varphi}-(q_L+2\lambda a_L)\dr_x\varphi\\
%    &+(g_L-\lambda q_L-\lambda\dr_xa_L+\lambda^2a_L)\varphi=\chi_{\lambda}[L,q,g]\varphi&\text{in $\Er\times\Theta$},\\
%    \nu\cdot\nabla_{\theta}\varphi&=0&\text{over $\Er\times\dr\Theta$},\\
%    \varphi&&\text{is $1$-periodic in $x$},\\
%    \varphi&>0&\text{in $\Er\times\Theta$}.
%  \end{aligned}
%  \right.
%\end{equation*}
%In particular, $k_L^{\lambda}=\chi_{\lambda}[L,0,r].$
%We let $\vect{q}(x)\in\Er^{1+P}$ be the vector corresponding to $q(x)$: \[\vect{q}(x):=\pmat{q(x)\\0\\\vdots\\0},\] %\\\vdots\\0},\]
%so that $\vect{q}\cdot\nabla\varphi^{\lambda}=q\dr_x\varphi^{\lambda}$ (here and below, we let $\nabla=\nabla_{(x,\theta)}$ be the gradient in all directions).
We let $A(x,\theta)$ be the total diffusion matrix in $\Er^{1+P}$, that is:
\[A(x,\theta)=
\pmat{a(x,\theta)&0&\cdots&0\\
  0&\mu(x,\theta)&\ddots&\vdots&\\
  \vdots&\ddots&\ddots&0\\
  0&\cdots&0&\mu(x,\theta)}.\]
The main tool of the proof of Proposition~\ref{ppn:LKL_nondecreasing} is the following theorem, which allows one to turn the study of a nonsymmetric operator into the study of a symmetric operator.

\begin{thm}[Theorem~2.2 in \cite{Nad10}]\label{thm:nadin}
  Let $q\in \mcc^{0,\alpha}(\Er\times\overline{\Theta})$ and $g\in \mcc^{0,\alpha}(\Er\times\overline{\Theta})$ be $1$-periodic in $x$.
  Let \[\mce:=\acc{\beta\in\mcc^{1,\alpha}(\Er\times\overline{\Theta})\tq\beta\text{ is $1$-periodic in $x$ and $\nu\cdot\nabla_{\theta}\beta=0$ on $\Er\times\dr\Theta$}}.\]
  Then:
  \begin{equation*}
  \chi[q,g]=\min_{\beta\in \mce}\chi\cro{0,\, g+{\nabla_{(x,\theta)}\beta\cdot A\nabla_{(x,\theta)}\beta}-q\dr_x\beta+\frac{\dr_xq}{2}}.
\end{equation*}
  
\end{thm}

The proof of Theorem~\ref{thm:nadin} is almost the same as that of \cite[Theorem~2.2]{Nad10}.
To make it work in our context, the Neumann boundary condition on $\Er\times\dr\Theta$, satisfied by the principal eigenfunction and by the elements of $\mce$, is central. 
Using Theorem~\ref{thm:nadin}, we may conclude in exactly the same way as in the proof of \cite[Proposition~4.1]{Nad10} that $L\mapsto k_L^{\lambda}$ is nondecreasing.

Thus Proposition~\ref{ppn:LKL_nondecreasing}, and with  Theorem~\ref{thm:freidlin_gartner},
we conclude that $L\mapsto c_L$ is nondecreasing.

\subsection{When the period goes to infinity (Theorem~\ref{thm:L_to_inf})}\label{ss:period_inf}

We let $H(x)=H_m(x)$ be the local principal eigenvalue, defined by the local eigenproblem
\begin{equation*}
  \left\{
  \begin{aligned}
    \nabla_{\theta}\cdot\pth{\mu\nabla_{\theta}\psi}+r(x,\cdot)\psi&=H(x)\psi&\text{in $\Theta$},\\
    \nu\cdot\nabla_{\theta}\psi&=0&\text{over $\dr\Theta$},\\
    \psi&>0&\text{in $\Theta$}.
  \end{aligned}
  \right.
\end{equation*}
The index $m$ being fixed to $1$, we drop it in this section for clarity.
We set $\displaystyle M:=\max_{x\in\Er}H(x)$ and we let $j:\ifo{M,+\infty}\to\Er$ be the function defined by
\[j(k)=\int_0^1\sqrt{\frac{k-H(x)}{a(x)}}\de x,\esp k\goq M.\]
Note that $j$ is a bijection from $\ifo{M,+\infty}$ to $\ifo{j(M),+\infty}$. We will prove the following proposition.

  \begin{ppn}\label{ppn:L_to_inf}
    Assume that $k^{0}_{L_0}>0$ for some $L_0>0$. Assume that $a$ is independent of $\theta$ and that $\mu$ is independent of $x$. Then $c_{L}$ converges as $L\to+\infty$ and
    \[\lim_{L\to +\infty}c_{L}=\inf_{\lambda\goq M}\frac{j^{-1}(\lambda)}{\lambda}.\]
  \end{ppn}
  
Theorem~\ref{thm:L_to_inf} is a consequence of the Proposition~\ref{ppn:L_to_inf} and of~\cite[Theorem~2.3]{HNR11}, which says that
\[\lim_{L\to+\infty}c^1(a_L,H_L)=\inf_{\lambda\goq M}\frac{j^{-1}(\lambda)}{\lambda}.\]
  To prove Proposition~\ref{ppn:L_to_inf}, we first prove the following lemma.
  
\begin{lem}\label{lem:use_qsd}
  Let $\lambda>0$.
  Let the assumptions of Proposition~\ref{ppn:L_to_inf} hold.
  Then $k_L^{\lambda}$ converges as $L\to+\infty$ and
  \begin{enumerate}[label=$(\roman*)$]
  \item If $\lambda\goq j(M)$, then
    \[\lim_{L\to+\infty}k_L^{\lambda}= j^{-1}(\lambda).\]
  \item If $\lambda<j(M)$, then 
  \[\lim_{L\to+\infty}k_L^{\lambda}= M.\]
  \end{enumerate}
%  If $\lambda\goq j(M)$, then $k_L^{\lambda}$ converges as $L\to+\infty$ and
%  \[\lim_{L\to+\infty}k_L^{\lambda}= j^{-1}(\lambda).\]
%  If $\lambda<j(M)$, then $k_L^{\lambda}$ converges as $L\to+\infty$ and
%  \[\lim_{L\to+\infty}k_L^{\lambda}= M.\]
\end{lem}
\begin{proof}
  Let $\varphi_L^{\lambda}$ be a positive principal eigenfunction associated to the Freidlin-Gärtner operator $\mcl_L^{\lambda}$.
  By considering the $L$-periodic version of $\varphi_L^{\lambda}$,
  which is $\varphi:(x,\theta)\mapsto\varphi_L^{\lambda}\pth{\frac{x}{L},\theta}$,
  we see that $k_L^{\lambda}$ is also a principal eigenvalue associated to the eigenproblem
\begin{equation*}
  \left\{
  \begin{aligned}
    \frac{1}{L^2}\dr_x(a\dr_x\varphi)+\nabla_{\theta}\cdot\pth{\mu\nabla_{\theta}\varphi}-\frac{2\lambda a}{L}\dr_{x}\varphi\esp&\\
    +\pth{r(x,\theta)+\lambda^2a-\lambda\frac{\dr_xa}{L}}\varphi&=k^\lambda_{L}\varphi&\text{in $\Er\times\Theta$},\\
    \nu\cdot\nabla_{\theta}\varphi&=0&\text{over $\Er\times\dr\Theta$},\\
    \varphi&>0&\text{in $\Er\times\Theta$},\\
    \varphi&&\text{is $1$-periodic in $x$}.
  \end{aligned}
  \right.
\end{equation*}
Therefore, $k^{\lambda}_L=\overline{k}^{\lambda}_L+o_{L\to+\infty}(1)$,
where $\overline{k}^{\lambda}_L$ is the principal eigenvalue associated to the same problem but without $-\lambda\frac{\dr_xa}{L}\varphi$:
\begin{equation*}
  \left\{
  \begin{aligned}
    \frac{1}{L^2}\dr_x(a\dr_x\varphi)+\nabla_{\theta}\cdot\pth{\mu\nabla_{\theta}\varphi}-\frac{2\lambda a}{L}\dr_{x}\varphi\esp&\\
    +\pth{r(x,\theta)+\lambda^2a}\varphi&=k^\lambda_{L}\varphi&\text{in $\Er\times\Theta$},\\
    \nu\cdot\nabla_{\theta}\varphi&=0&\text{over $\Er\times\dr\Theta$},\\
    \varphi&>0&\text{in $\Er\times\Theta$},\\
    \varphi&&\text{is $1$-periodic in $x$}.
  \end{aligned}
  \right.
\end{equation*}
We now apply~\cite[Theorem 1.3]{b23}. To recover the notations there, we let~$\eps=1/L$, and replace the slow variable by~$x$ (instead of~$y$) and the fast variable by~$\theta$ (instead of~$z$).
We rewrite $j$ as
\[j(k)=\int_0^1\sqrt{\frac{k-(H(x)+\lambda^2a(x))+\frac{(2\lambda a(x))^2}{4a(x)}}{a(x)}}\de x.\]
Now,~\cite[Theorem 1.3]{b23} says that as $L\to+\infty$ (\ie $\eps\to 0$),
\[\overline{k}_L^{\lambda}\to J^{-1}\pth{\abs{\int_0^1\frac{2\lambda a(x)}{2a(x)}}}=j^{-1}(\lambda)\]
if $\lambda\goq J(M)=j(M)$, and 
\[\overline{k}_L^{\lambda}\to M\]
otherwise.
Hence the same convergence holds for $k^{\lambda}_L$.
\end{proof}

We are now ready to conclude the proof of Proposition~\ref{ppn:L_to_inf}. 
\begin{proof}[Proof of Proposition~\ref{ppn:L_to_inf}]
  By the persistence assumption that $k^{0}_{L_0}>0$, and by Proposition~\ref{ppn:LKL_nondecreasing}, there exists $\delta>0$ such that for each $L\goq L_0$: $k^0_L>\delta$.
  This and the lower bound in Lemma~\ref{lem:estimate_klambda} imply that for each $L\goq L_0$, there exists $\lambda_L>0$ such that 
  \[\inf_{\lambda>0}\frac{k_L^{\lambda}}{\lambda}=\frac{k_L(\lambda_L)}{\lambda_L},\]
  and that there exist $\lambda_->0$ and $\lambda_+>\lambda_-$ such that for all $L>0$, $\lambda_-\loq\lambda_L\loq \lambda_+$.
  Thus, there exist an increasing sequence $L_n\to+\infty$ and $\lambda_{\infty}>0$ such that, as $n\to+\infty$, $\lambda_{L_n}\to\lambda_{\infty}$.
  Since the sequence of functions $\lambda\mapsto k_{L_n}^{\lambda}$ converges increasingly, as $n\to+\infty$, to its limit,
  the convergence is locally uniform.
  We then have,
  \[\liminf_{n\to+\infty}c_{L_n}
  =\liminf_{n\to+\infty}\frac{k_{L_n}(\lambda_{L_n})}{\lambda_{L_n}}
  =\liminf_{n\to+\infty}\frac{k_{L_n}(\lambda_{\infty})}{\lambda_{\infty}}.
  \]
  Therefore,
  \[\liminf_{n\to+\infty}c_{L_n}\goq\inf_{\lambda>0}\lim_{L\to+\infty}\frac{k^{\lambda}_L}{\lambda}.\]
  Finally, for all $n\goq 1$, and for all $\lambda>0$, $c_{L_n}\loq\frac{k_{L_n}^{\lambda}}{\lambda}$.
  Taking the limit $n\to+\infty$, we obtain, for all $\lambda>0$,
  \[\limsup_{n\to+\infty}c_{L_n}\loq\lim_{n\to+\infty}\frac{k_{L_n}^{\lambda}}{\lambda}.\]
  With the above, we obtain that $c_{L_n}$ converges for all sequence $L_n\to+\infty$.
  Hence $c_L$ also converges and:
  \[\lim_{L\to+\infty}c_{L}=\inf_{\lambda>0}\lim_{n\to+\infty}\frac{k_{L_n}^{\lambda}}{\lambda}.\]
  To conclude, by Lemma~\ref{lem:use_qsd}, for $\lambda\goq j(M)$, we have
  \[\lim_{L\to+\infty}k_{L}^{\lambda}=j^{-1}(\lambda),\]
  and for $\lambda<j(M)$, we have
  \[\lim_{L\to+\infty}k_{L}^{\lambda}\goq M=\lim_{L\to+\infty}k_L^{j(M)}.\]
  Thus
  \[\lim_{L\to+\infty}c_{L}=\inf_{\lambda>0}\lim_{L\to+\infty}\frac{k_{L}^{\lambda}}{\lambda}=\inf_{\lambda>j(M)}\lim_{L\to+\infty}\frac{k_{L}^{\lambda}}{\lambda}.\]
  This concludes the proof.
\end{proof}

\begin{proof}[Proof of Equation~\eqref{eq:cH_vs_cinf}]
  Now, let us prove that~\eqref{eq:cH_vs_cinf} holds when $a$ is a constant.
  We assume without loss of generality that $a=1$. 
  On the one hand, we have by Proposition~\ref{ppn:L_to_inf},
  \begin{align*}
    \lim_{L\to+\infty}c_L&=\inf_{\lambda\goq j(M)}\frac{j^{-1}(\lambda)}{\lambda}=\inf_{k\goq M}\frac{k}{j(k)}\\
    &=\inf_{k\goq M}k\pth{\int_0^1\sqrt{k-H(x)}\de x}^{-1}.
  \end{align*}
  By Jensen's inequality, we obtain:
   \begin{align*}
     \lim_{L\to+\infty}c_L&\goq\inf_{k\goq M}\sqrt{k^2\pth{\int_0^1(k-H(x))\de x}^{-1}}\\
     &=\sqrt{\inf_{k\goq M}\pth{\frac{k^2}{k-\widetilde{H}^A}}},
   \end{align*}
   where $\widetilde{H}^A=\int_0^1H(x)\de x$ is the arithmetic mean of $H$.
   But:
   \begin{align*}
     \inf_{k\goq M}\pth{\frac{k^2}{k-\widetilde{H}^A}}&=
     \inf_{k\goq M}\pth{2\widetilde{H}^A+(k-\widetilde{H}^A)+\frac{(\widetilde{H}^A)^2}{k-\widetilde{H}^A}}\\
     &=2\widetilde{H}^A+\inf_{k\goq M-\widetilde{H}^A}\pth{k+\frac{(\widetilde{H}^A)^2}{k}}
     \goq 4\widetilde{H}^A.
   \end{align*}
   Hence,
   \begin{align*}
     \lim_{L\to+\infty}c_L&\goq 2\sqrt{\widetilde{H}^A}.
   \end{align*}
   On the other hand, applying Jensen's inequality twice,
   \begin{align*}
     c_H&=2\pth{\int_0^1\frac{\de x}{\sqrt{H(x)}}}^{-1}\loq 2\int_0^1\sqrt{H(x)}\de x\loq 2\sqrt{\int_0^1H(x)\de x}=2\sqrt{\widetilde{H}^A}.
   \end{align*}
   The inequality is strict if $H$ is not constant. 
   The inequality~\eqref{eq:cH_vs_cinf} is proved.

\end{proof}

\subsection{When the period goes to zero (Theorem~\ref{thm:L_to_0})}\label{ss:period_0} 

In order to compute the limit of $c_L$ as $L\to 0$, we first compute the limit as $L\to 0$ of $k_L^{\lambda}$.
As in the introduction, we define the \emph{homogenised principal eigenvalue $\widetilde{k}^{\lambda}$} and an associated homogenised principal eigenfunction $\psi^{\lambda}\in\mcc^{2}(\overline{\Theta})$ by 
\begin{equation*}%\label{eq:homo_eig_problem}
  \left\{
  \begin{aligned}
    \nabla_{\theta}\cdot\pth{\widetilde{\mu}^A\nabla_{\theta}\psi}+\pth{\widetilde{r}^A+\lambda^2\widetilde{a}^H}\psi&=\widetilde{k}^{\lambda}\psi&\comment{in $\Theta$},\\
    \nu\cdot\nabla_{\theta}\psi&=0&\text{over $\dr\Theta$},\\
    \psi&>0&\comment{in $\Theta$}.
  \end{aligned}
  \right.
\end{equation*}

\begin{lem}\label{lem:k_Lto0}
  Fix $\lambda\in\Er$. Then, as $L\to 0$,
  \[k_L^{\lambda}\to\widetilde{k}^{\lambda}.\]
\end{lem}

\begin{proof}
  Since $\lambda$ is fixed, we shall omit the superindex $\lambda$ in the notations. For all $L>0$, we let $\varphi_L:=\varphi_L^{\lambda}>0$ and $k_L:=k_L^{\lambda}$ solve
  \begin{equation}\label{eq:eig_pb_Lto0}
  \left\{
  \begin{aligned}
    \dr_x\pth{a_L\dr_x\varphi_L}+\nabla_{\theta}\cdot\pth{\mu_L\nabla_{\theta}\varphi_L}-2\lambda a_L\dr_x\varphi_L\esp&\\
    +\pth{r_L+\lambda^2a_L-\lambda\dr_xa_L}\varphi_L&=k_L\varphi_L&\text{in $\Theta$},\\
    \nu\cdot\nabla_{\theta}{\varphi_L}&=0&\text{over $\dr\Theta$},\\
    {\varphi_L}&&\text{is $L$-periodic in $x$},
  \end{aligned}
  \right.
  \end{equation}
  with normalisation
  \[\int_{[0,1]\times\Theta}\varphi_L^2=1.\]
  (Note that we integrate on $[0,1]\times\Theta$, so for $L\neq 1$ we do not integrate exactly on a periodic cell.)
  We let $\Phi_L(x,\theta):=e^{-\lambda x}\varphi_L(x,\theta)$, so that $\Phi_L$ solves
  \begin{equation}\label{eq:equation_Phi_Lto0}
  \left\{
  \begin{aligned}
    \dr_x\pth{a_L\dr_x\Phi_L}+\nabla_{\theta}\cdot\pth{\mu_L\nabla_{\theta}\Phi_L}+r_L\Phi_L&=k_L^{\lambda}\Phi_L&\text{in $\Theta$},\\
    \nu\cdot\nabla_{\theta}{\Phi_L}&=0&\text{over $\dr\Theta$}.
  \end{aligned}
  \right.
\end{equation}

  \paragraph{Step 1. Convergence of $\varphi_L$ to a function $\varphi_{(0)}$ independent of $x$.} 
  %We note ${\omega=\cro{0,1}\times \Theta}$ and $\omega_L=\cro{0,L\times\ceil{1/L}}\times\Theta$.
  First of all, we show that $(\varphi_L)_L$ is bounded in $W_2^1([0,1]\times\Theta)$.
  For $L>0$, we note $\widehat{L}:=L\times\ceil{1/L}$. We point out that $\widehat{L}\to 1$ as $L\to 0$.
  We multiply the first line of \eqref{eq:eig_pb_Lto0} by~$\varphi_L$, we integrate on ${[0,\widehat{L}]\times\Theta}$ and we integrate by parts. Using the $L$-periodicity in $x$ and $\nu\cdot\nabla_{\theta}\varphi_L=0$ on $\Er\times\dr\Theta$, we get
  \begin{multline*}
    -\int_{{[0,\widehat{L}]\times\Theta}} a_L\abs{\dr_x\varphi_L}^2-\int_{{[0,\widehat{L}]\times\Theta}}\mu_L\abs{\nabla_{\theta}\varphi_L}^2-2\lambda\int_{{[0,\widehat{L}]\times\Theta}} a_L\varphi_L\dr_x\varphi_L
    \\+\int_{{[0,\widehat{L}]\times\Theta}}\pth{r_L+\lambda^2a_L-\lambda\dr_xa_L}\varphi_L^2
    =k_L\int_{{[0,\widehat{L}]\times\Theta}} \varphi_L^2.    
  \end{multline*}
  Thus, replacing ${[0,\widehat{L}]\times\Theta}$ by $[0,1]\times\Theta$,
  \begin{multline*}
    -\int_{[0,1]\times\Theta} a_L\abs{\dr_x\varphi_L}^2-\int_{[0,1]\times\Theta}\mu_L\abs{\nabla_{\theta}\varphi_L}^2-2\lambda\int_{[0,1]\times\Theta} a_L\varphi_L\dr_x\varphi_L\\
    +\int_{[0,1]\times\Theta}\pth{r_L+\lambda^2a_L-\lambda\dr_xa_L}\varphi_L^2
    =k_L\int_{[0,1]\times\Theta} \varphi_L^2+o_{L\to 0}(1).    
  \end{multline*}
 Recall the assumption that there exists a constant $\eta>0$ such that $a_L\goq\eta$ and $\mu_L\goq\eta$.
 We get, rearranging the terms:
 \begin{multline*}
   \eta\int_{[0,1]\times\Theta}\abs{\nabla\varphi_L}^2\loq
   -\lambda\int_{[0,1]\times\Theta} \varphi_L^2\dr_xa_L\\
   -2\lambda\int_{[0,1]\times\Theta} a_L\varphi_L\dr_x\varphi_L+\int_{[0,1]\times\Theta}\pth{r_L+\lambda^2a_L}\varphi_L^2-k_L\int_{[0,1]\times\Theta} \varphi_L^2+o_{L\to 0}(1),   
 \end{multline*}
  with $\nabla\varphi_L=(\dr_x\varphi_L,\nabla_{\theta}\varphi_L)$.
  We have:
  \begin{align*}
    -2\int_{{[0,\widehat{L}]\times\Theta}} a_L\varphi_L\dr_x\varphi_L-\int_{{[0,\widehat{L}]\times\Theta}} \varphi_L^2\dr_xa_L&=-\int_{{[0,\widehat{L}]\times\Theta}} a_L\dr_x(\varphi_L^2)-\int_{{[0,\widehat{L}]\times\Theta}} \varphi_L^2\dr_xa_L\\
    &=\int_{{[0,\widehat{L}]\times\Theta}} \varphi_L^2\dr_xa_L-\int_{{[0,\widehat{L}]\times\Theta}} \varphi_L^2\dr_xa_L\\
    &=0.
  \end{align*}
  We conclude that
  \[
    \int_{[0,1]\times\Theta}\abs{\nabla\varphi_L}^2
    \loq
    \frac{1}{\eta}\pth{(\dabs{r_L}_{\infty}+\lambda^2\dabs{a_L}_{\infty})\int_{[0,1]\times\Theta}\varphi_L^2-k_L\int_{[0,1]\times\Theta} \varphi_L^2}+o_{L\to 0}(1).
  \]
  Since $(k_L)_{L>0}$ is bounded (see Lemma~\ref{lem:estimate_klambda}) and $\dabs{\varphi_L}_{L^2([0,1]\times\Theta)}=1$, the family $(\varphi_L)_L$ is bounded in $W_2^1([0,1]\times\Theta)$.
%  Thus, there is a sequence $L_n\to 0$ such that $(\varphi_{L_n})$ converges in $L^2([0,1]\times\Theta)$ to a function $\varphi_{(0)}\in L^2([0,1]\times\Theta)$.
  Thus, by a compact injection theorem (here compactness is due to periodicity),
  there is a sequence $L_n\to 0$ such that $(\varphi_{L_n})$ converges in $L^2([0,1]\times\Theta)$ to a function $\varphi_{(0)}\in L^2([0,1]\times\Theta)$.
  For convenience, we note $\varphi_n:=\varphi_{L_n}$, $\Phi_n:=\Phi_{L_n}$ and $k_n:=k_{L_n}$.
  Since $\varphi_{n}$ is $L_n$-periodic in $x$, the limit $\varphi_{(0)}(x,\theta)=\varphi_{(0)}(\theta)$ must be independent of~$x$.
  We conclude that ${\Phi}_n$ converges in $L^2_{loc}(\Er\times\overline{\Theta})$ to the function 
\[\Phi_{(0)}(x,\theta):=e^{-\lambda x}\varphi_{(0)}(\theta).\]

  \paragraph{Step 2. Apply the oscillating test function method.} Since $(k_L)_{L>0}$ is bounded (see Lemma~\ref{lem:estimate_klambda}), there exists $k_{(0)}\in\Er$ such that, up to extraction of a subsequence of $(L_n)$,
  \[\lim_{n\to+\infty}k_n=k_{(0)}.\]
  Now, we show that in the classical sense:
  \begin{equation*}
    \widetilde{a}^H(\theta)\dr_{xx}\Phi_{(0)}+\nabla_{\theta}\cdot\pth{\widetilde{\mu}^A(\theta)\nabla_{\theta}\Phi_{(0)}}+\widetilde{r}^A(\theta)\Phi_{(0)}=k_{(0)}\Phi_{(0)},
  \end{equation*}
  and that $\nu\cdot\nabla_{\theta}\Phi_{(0)}=0$ on $\Er\times\dr\Theta$.
  As a first step,
  we let $\psi\in W^1_2(\Er\times{\Theta})$ be a test function such that $e^{-\lambda x}\psi$ is $1$-periodic in~$x$
  (so that $\psi\Phi_{(0)}$, $\psi\Phi_L$, $\psi\dr_x\Phi_{(0)}$ and $\psi\dr_x\Phi_K$ are $1$-periodic in $x$), and we show that
  \begin{equation}\label{eq:weak_formulation_psi}
    \int_{[0,1]\times\Theta}\widetilde{a}^H (\dr_x\Phi_{(0)})(\dr_x\psi)-\int_{[0,1]\times\Theta}\widetilde{\mu}^A\nabla_{\theta}\Phi_{(0)}\cdot\nabla_{\theta}\psi+\int_{[0,1]\times\Theta}\pth{\widetilde{r}^A-k_{(0)}}\Phi_{(0)}\psi=0.
  \end{equation}
  To this aim, we apply the oscillating test function method, which is classical in homogenisation theory~\cite{All12}. We let
  \[w(x,\theta):=\int_0^x\pth{\frac{\widetilde{a}^H(\theta)}{a(x',\theta)}-1}\de x'.\]
  We note that $w$ is $1$-periodic in $x$ and solves
  \begin{equation}
    \label{eq:cell_problem}\dr_x\cro{a(x,\theta)(1+\dr_xw(x,\theta))}=0.
  \end{equation}
  We set, for $n\goq 0$:
  \[\psi_n(x,\theta):=\psi(x,\theta)+L_n\dr_x\psi(x,\theta)w\pth{\frac{x}{L_n},\theta}.\]
  Since $\nu\cdot\nabla_{\theta}\Phi_n=0$ on $\Er\times\dr\Theta$, multiplying \eqref{eq:equation_Phi_Lto0} by $\psi_n$ and integrating by parts on $[0,1]\times\Theta$ gives: 
  \begin{multline}\label{eq:weak_formulation_psin}
    -\int_{[0,1]\times\Theta} a\pth{\frac{x}{L_n},\theta}(\dr_x\Phi_n)(\dr_x\psi_n)-\int_{[0,1]\times\Theta}\mu\pth{\frac{x}{L_n},\theta}\nabla_{\theta}\Phi_{n}\cdot\nabla_{\theta}\psi\\
    +\int_{[0,1]\times\Theta}\pth{r\pth{\frac{x}{L_n},\theta}-k_n}\Phi_n\psi_n
    =o_{L\to 0}(1).
  \end{multline}
  As above, the $o_{L\to 0}(1)$ term comes from the fact that the functions which we integrate by parts are  $L_n$-periodic in $x$ and may not be exactly $1$-periodic in $x$.

  Now, let us compute the limit as $n\to+\infty$ for each term of Equation \eqref{eq:weak_formulation_psin}. First, 
  \begin{align*}
    &\int_{[0,1]\times\Theta} a\pth{\frac{x}{L_n},\theta}(\dr_x\Phi_n)(\dr_x\psi_n)\nonumber\\
    &\esp=\int_{[0,1]\times\Theta} a\pth{\frac{x}{L_n},\theta}(\dr_x\Phi_n)(\dr_x\psi)\pth{1+\dr_xw\pth{\frac{x}{L_n},\theta}}\\
    &\esp\esp\esp+L_n\int_{[0,1]\times\Theta} a\pth{\frac{x}{L_n},\theta}(\dr_x\Phi_n)(\dr_{xx}\psi)w\pth{\frac{x}{L_n},\theta}.
  \end{align*}
  On the one hand, %since $w$ satisfies \eqref{eq:cell_problem}, there exists a constant $a^*>0$ such that:
%  \[{a\pth{\frac{x}{L_n},\theta}\pth{1+\dr_xw\pth{\frac{x}{L_n},\theta}}}\equiv a^*.\]
%  We obtain:
  \begin{equation*}
    \int_{[0,1]\times\Theta} a\pth{\frac{x}{L_n},\theta}(\dr_x\Phi_n)(\dr_x\psi)\pth{1+\dr_xw\pth{\frac{x}{L_n},\theta}}=\widetilde{a}^H\int_{[0,1]\times\Theta}\dr_x\Phi_n\dr_x\psi.
  \end{equation*}
  On the other hand, since $\Phi^L$ is bounded in $W_2^1([0,1]\times\Theta)$, the family
  \[\pth{\int_{[0,1]\times\Theta} a\pth{\frac{x}{L},\theta}(\dr_x\Phi^L)(\dr_{xx}\psi)w\pth{\frac{x}{L},\theta}}_{0<L<1}\]
  is bounded. Thus, 
  \begin{equation}\label{eq:lim_L0_a}
    \lim_{n\to+\infty}\int_{[0,1]\times\Theta} a\pth{\frac{x}{L_n},\theta}(\dr_x\Phi_n)(\dr_x\psi_n)=\widetilde{a}^H\int_{[0,1]\times\Theta}\dr_x\Phi_{(0)}\dr_x\psi.
  \end{equation}
  Second, an integration by part gives:
  \[\int_{[0,1]\times\Theta} \mu\pth{\frac{x}{L_n},\theta}\nabla_{\theta}\Phi_{n}\cdot\nabla_{\theta}\psi=-\int_{[0,1]\times\Theta}\Phi_n\nabla_{\theta}\cdot\pth{ \mu\pth{\frac{x}{L_n},\theta}\nabla_{\theta}\psi_n}.\]
  As $n\to +\infty$, the right-hand side goes to
\[  -\int_{[0,1]\times\Theta}\Phi_{(0)}\nabla_{\theta}\cdot\pth{ \widetilde{\mu}^A\nabla_{\theta}\psi}.\]
  Thus, by another integration by part:
  \begin{equation}\label{eq:lim_L0_mu}
        \lim_{n\to +\infty}\int_{[0,1]\times\Theta} \mu\pth{\frac{x}{L_n},\theta}\nabla_{\theta}\Phi_{n}\cdot\nabla_{\theta}\psi=\int_{[0,1]\times\Theta}\widetilde{\mu}^A\pth{\theta}(\nabla_{\theta}\Phi_{(0)})(\nabla_{\theta}\psi).
  \end{equation}
  Finally,
  \begin{equation}\label{eq:lim_L0_r}
    \lim_{n\to+\infty}\int_{[0,1]\times\Theta} \pth{r\pth{\frac{x}{L_n},\theta}-k_{(0)}}\Phi_n\psi_n=\int_{[0,1]\times\Theta}\pth{\widetilde{r}^A(\theta)-k_{(0)}}\Phi_{(0)}\psi.
  \end{equation}
  Using Equations \eqref{eq:weak_formulation_psin}, \eqref{eq:lim_L0_a}, \eqref{eq:lim_L0_mu} and \eqref{eq:lim_L0_r},
  we conclude that the weak formulation \eqref{eq:weak_formulation_psi} holds
  for all test functions $\psi\in W^1_2(\Er\times\Theta)$
  such that $e^{-\lambda x}\psi$ is $1$-periodic in~$x$.
  Hence, by the notes at the end of \cite[Chapter~8]{GilTru01},
  the function~$\Phi_{(0)}$ is the unique weak solution of
    \sys{\label{eq:psi_bar_homog}
      \widetilde{a}^H\dr_{xx}\Phi_{(0)}+\nabla_{\theta}\cdot\pth{\widetilde{\mu}^A\nabla_{\theta}\Phi_{(0)}}+\widetilde{r}^A\Phi_{(0)}&=k_{(0)}\Phi_{(0)}&\text{in $\Er\times\Theta$},\\
      \nu\cdot\nabla_{\theta}\Phi_{(0)}&=0&\text{over $\Er\times\dr\Theta$}.
      }

    \paragraph{Step 3. Conclusion.}
    By the regularity of the coefficients of \eqref{eq:psi_bar_homog} and the uniqueness of the classical solution, we conclude that $\Phi_{(0)}\in\mcc^{2}(\Er\times\overline{\Theta})$, so \eqref{eq:psi_bar_homog} is solved in the classical sense.
    We have also $\Phi_{(0)}(x,\theta)=e^{-\lambda x}\varphi_{(0)}(\theta)$.
  Therefore $\varphi_{(0)}\in\mcc^{2}(\Er\times\overline{\Theta})$ solves: 
  \[\nabla_{\theta}\cdot\pth{\mu\nabla_{\theta}\varphi_{(0)}}+\pth{\lambda^2\widetilde{a}^H+\widetilde{r}^A}\varphi_{(0)}=k_{(0)}\varphi_{(0)},\]
  with $\nu\cdot\nabla_{\theta}\varphi_{(0)}=\nu\cdot\nabla_{\theta}\Phi_{(0)}=0$ on $\Er\times\dr\Theta$.
  Finally, since $\varphi_L\to\varphi_{(0)}$ in $L^2(\Er\times\Theta)$, we have $\varphi_{(0)}\goq 0$
  and $\varphi_{(0)}\not\equiv 0$.
  Using the strong comparison principle, we obtain: $\varphi_{(0)}> 0$.
  
  Therefore, by the uniqueness of a principal eigenvalue, we have $k_{(0)}=\widetilde{k}$.
  To conclude, the family $(k_L)_{L>0}$ is bounded (see Lemma~\ref{lem:estimate_klambda}) and $\widetilde{k}$ is the only possible limit of a subsequence as $L\to 0$; hence the whole family converges to $\widetilde{k}$ as $L\to 0$.

\end{proof}

From Lemma~\ref{lem:k_Lto0}, we show the first equality in Theorem~\ref{thm:L_to_0}. 

\begin{proof}[Proof of Theorem~\ref{thm:L_to_0}, part 1]
  By Lemma~\ref{lem:k_Lto0} and the assumption that $\widetilde{k}^0>0$, there exist $\delta>0$ and $L_0>0$ such that for all $L\in\ioo{0,L_0}$,
  \[k^0_{L}>\delta.\]
  Thus the function $L\mapsto c_L$ is nondecreasing on $\ioo{0,+\infty}$ (see Proposition~\ref{ppn:inter_L}). Thus
  \[\lim_{L\to 0}c_L=\inf_{L>0}c_L.\]
  Thus:
  \[\lim_{L\to 0}c_L=\inf_{L>0}\pth{\inf_{\lambda>0}\frac{k_L^{\lambda}}{\lambda}}=\inf_{\lambda>0}\pth{\inf_{L>0}\frac{k_L^{\lambda}}{\lambda}}.\]
  Since, for each $\lambda>0$, the function $L\mapsto k_L^{\lambda}$ is nondecreasing (see Proposition~\ref{ppn:LKL_nondecreasing}), we have:
  \[\inf_{L>0}\frac{k_L^{\lambda}}{\lambda}=\frac{\widetilde{k}^{\lambda}}{\lambda}.\]
  Thus:
  \[\lim_{L\to 0}c_L=\inf_{\lambda>0}\frac{\widetilde{k}^{\lambda}}{\lambda}.\]
\end{proof}

Finally, we show the second equality of Theorem~\ref{thm:L_to_0} by using the Rayleigh formula.

\begin{proof}[Proof of Theorem~\ref{thm:L_to_0}, part 2]
  For $\psi\in H^1(\Theta)$, we let
  \[A_{\psi}=\int_{\Theta}{\widetilde{a}^H\psi^2},\esp R_{\psi}=\int_{\Theta}\pth{\widetilde{r}^A\psi^2-\widetilde{\mu}^A\abs{\nabla\psi}^2}.\]
  Let us note
  \[c_{var}:=\sup_{\psi\in H^1(\Theta),\,\|\psi\|_{L^2(\Theta)}=1,\,R_{\psi}>0}2\sqrt{A_{\psi}R_{\psi}}\]
  and
  \[c_0:=\inf_{\lambda>0}\frac{\widetilde{k}^{\lambda}}{\lambda}.\]
  We will show that $c_0=c_{var}$. 
  
  For $\lambda> 0$, let $\widetilde{k}^{\lambda}$ be the homogenised principal eigenvalue (defined by \eqref{eq:eig_pb_Lto0}). Recall the assumption that for all $L\in\ioo{0,L_0}$, $k^0_{m,L}>\delta$.
  By Lemma~\ref{lem:k_Lto0}, therefore, we must have $\widetilde{k}^0>0$;
  due to the Rayleigh formula for $\widetilde{k}^0$, there exists $\psi$ such that $R_{\psi}>0$. 
  Therefore, the set on which the $\sup$ is taken in the definition of $c_{var}$ is nonempty, so $c_{var}$ is well-defined.
  
By the Rayleigh formula, we have
  \[\widetilde{k}^{\lambda}=\max_{{\psi\in H^1(\Theta),\,\|\psi\|_{L^2(\Theta)}=1}}\pth{R_{\psi}+\lambda^2 A_{\psi}}.\]
  Hence
  \begin{equation}\label{eq:ex_c0*}
    c_0=\inf_{\lambda>0}\cro{\max_{\psi\in H^1(\Theta),\,\|\psi\|_{L^2(\Theta)}=1}\pth{\frac{1}{\lambda}R_{\psi}+\lambda A_{\psi}}}.
  \end{equation}
  Moreover, for all $\psi\in H^1(\Theta)$,
  \[\min_{\lambda>0}\pth{\frac{1}{\lambda}R_{\psi}+\lambda A_{\psi}}= 2\sqrt{A_{\psi}R_{\psi}}>0\]
  if $R_{\psi}>0$,
  and
  \[\min_{\lambda>0}\pth{\frac{1}{\lambda}R_{\psi}+\lambda A_{\psi}}= -\infty\]
  if $R_{\psi}<0$.
  Thus
  \begin{equation}\label{eq:ex_c'}
    c_{var}=\sup_{\psi\in H^1(\Theta),\,\|\psi\|_{L^2(\Theta)}=1,\,R_{\psi}>0}\cro{\min_{\lambda>0}\pth{\frac{1}{\lambda}R_{\psi}+\lambda A_{\psi}}}.
  \end{equation}
  First, we show that $c_{var}\loq c_0$.
  For all $\psi_0$ such that $\|\psi\|_{L^2(\Theta)}=1$ and $R_{\psi_0}>0$,
  and for all $\lambda_0>0$, we have:
  \[\min_{\lambda>0}\pth{\frac{1}{\lambda}R_{\psi_0}+\lambda A_{\psi_0}}\loq \max_{\psi\in H^1(\Theta),\,\|\psi\|_{L^2(\Theta)}=1,\,R_{\psi}>0}\pth{\frac{1}{\lambda_0}R_{\psi}+\lambda_0 A_{\psi}}.\]
  By \eqref{eq:ex_c0*} and \eqref{eq:ex_c'}, we obtain: $c_{var}\loq c_0$.
  Second, we show that $c_{var}\goq c_0$. For $\psi\in H^1(\Theta)$ such that $R_{\psi}>0$, let $\lambda_{\psi}$ be the argument that reaches the minimum of
  \[\lambda\mapsto\frac{R_{\psi}}{\lambda}+\lambda A_{\psi}.\]
  We then have:
  \begin{align*}
    c_0&\loq\sup_{\psi\in H^1(\Theta),\,\|\psi\|_{L^2(\Theta)}=1,\,R_{\psi}>0}\pth{\frac{R_{\psi}}{\lambda_{\psi}}+\lambda_{\psi} A_{\psi}}\\
  &=\sup_{\psi\in H^1(\Theta),\,\|\psi\|_{L^2(\Theta)}=1,\,R_{\psi}>0}\min_{\lambda>0}\pth{\frac{R_{\psi}}{\lambda}+\lambda A_{\psi}}=c_{var}.    
  \end{align*}
  Hence $c_{var}=c_0$, which concludes.
\end{proof}

\section{Behaviour of the spreading speed with respect to the mutation coefficient}\label{s:mutation}

This part is devoted to the proofs of Theorem~\ref{thm:mutation} (limits of $c_{m,L}$ as $m\to 0$ and $m\to+\infty$), of Theorem~\ref{thm:bio_csq} (optimal mutation coefficient), and of Propositions~\ref{ppn:limit_dec_mu} and~\ref{ppn:persistence_decr}, which give various properties of the dependence in the parameter $m$ of the spreading speed.

We shall use specific notations for principal eigenvalues for operators without phenotypic variable. 
If $A\in\mcc^{1,\alpha}(\Er)$ and $R\in\mcc^{0,\alpha}(\Er)$ are two $1$-periodic functions, we let $\kappa^1(\lambda;A,R)$
be the principal eigenvalue corresponding to the operator
\[\mcl^{1,\lambda}\phi:=\dr_x(A\dr_x\phi)-2\lambda A\dr_x\phi+\pth{R+\lambda^2A-\lambda\dr_xA}\phi\]
acting on $1$-periodic functions.
That is, $\kappa^1(\lambda;A,R)$ is the only possible value such that there exists a $1$-periodic function $\phi\in\mcc^2(\Er)$ satisfying
\[\mcl^{1,\lambda}\phi=\kappa^1(\lambda;A,R)\phi,\esp\phi>0.\]
Finally, for a function $f:\Er\times\Theta\to\Er$, we let
\[\overline{f}(x):=\frac{1}{\abs{\Theta}}\int_{\Theta}f(x,\theta)\de\theta\]
be the arithmetic mean of $f$ in the $\theta$-variable.

In Subsection~\ref{ss:m_to_inf}, we compute the limit of $c_m$ as $m\to+\infty$. In Subsection~\ref{ss:m_to_0}, we compute the limit of $c_m$ as $m\to 0$ and conclude the proof of Theorem~\ref{thm:mutation}. In Subsection~\ref{ss:other}, we prove Propositions~\ref{ppn:limit_dec_mu} and~\ref{ppn:persistence_decr}. Finally, in Subsection~\ref{ss:proof_bio}, we prove Theorem~\ref{thm:bio_csq}. 

\subsection{When the mutation coefficient goes to infinity (Theorem~\ref{thm:mutation}, part~$(i)$)}\label{ss:m_to_inf}

For notational convenience, we shall assume throughout Subsections~\ref{ss:m_to_inf} and~\ref{ss:m_to_0} that $L=1$. 
The notations thus simplify:
\[c_m=c_{m,L},\esp k^{\lambda}_m=k^{\lambda}_{m,L}.\]
In order to compute the limit of $c_m$ as $m\to +\infty$, we first compute the limit of $k_m^{\lambda}$ as $m\to +\infty$.

\begin{lem}\label{lem:limits_k_muinfty}
  Let $\lambda\in\Er$. Under the assumptions of Theorem~\ref{thm:mutation} (part~$(i)$), as $m\to +\infty$,
  \[k_m^{\lambda}\to \kappa^1\pth{\lambda;{a},\overline{r}}.\]
\end{lem}

\begin{proof}
  We assume that $a=a(x)$ and $\mu=\mu(x)$ are independent of $\theta$. 

  For all $m>0$, we let $\varphi_m>0$ be the principal eigenfunction associated to $k_m^{\lambda}$ such that $\dabs{\varphi_m}_{L^2([0,1]\times\Theta)}=1$, \ie  $\varphi_m$ is positive, $1$-periodic in $x$, and satisfies
  \begin{equation}\label{eq:varphim_mtoinf}
  \left\{
  \begin{aligned}
    \dr_x(a\dr_x\varphi_m)+m\nabla_{\theta}\cdot\pth{\mu\nabla_{\theta}\varphi_m}-2\lambda a\dr_x\varphi
    +(r+\lambda^2a-\lambda\dr_xa)\varphi_m&=k_m^{\lambda}\varphi_m&\text{in $\Er\times\Theta$},\\
    \nu\cdot\nabla_{\theta}\varphi_m&=0&\text{over $\Er\times\dr\Theta$}.
  \end{aligned}
  \right.
  \end{equation}
  First, we show that $(\varphi_m)_{m>1}$ is bounded in $W_2^1([0,1]\times\Theta)$. Multiplying the first line of \eqref{eq:varphim_mtoinf} by~$\varphi_m$, integrating on $[0,1]\times\Theta$, and integrating by parts, yields (using that the functions are $1$-periodic in~$x$ and $\nu\cdot\nabla_{\theta}\varphi_m=0$ on $\Er\times\dr\Theta$):
  \begin{multline*}
    -\int_{[0,1]\times\Theta} a\abs{\dr_x\varphi_m}^2-m\int_{[0,1]\times\Theta}\mu\abs{\nabla_{\theta}\varphi_m}^2-2\lambda\int_{[0,1]\times\Theta} a\varphi_m\dr_x\varphi_m\\
    +\int_{[0,1]\times\Theta}\pth{r+\lambda^2a-\lambda\dr_xa}\varphi_m^2=k_m\int_{[0,1]\times\Theta} \varphi_m^2.
  \end{multline*}
 Recall the assumption that there exists a constant $\eta>0$ such that $a\goq\eta$ and $\mu\goq\eta$.
 We get, rearranging the terms:
 \begin{multline*}
  \eta\pth{\int_{[0,1]\times\Theta}\abs{\dr_x\varphi_m}^2+m\int_{[0,1]\times\Theta}\abs{\nabla_{\theta}\varphi_m}^2}\\
  \loq -\lambda\int_{[0,1]\times\Theta} \varphi_m^2\dr_xa-2\lambda\int_{[0,1]\times\Theta} a\varphi_m\dr_x\varphi_m+\int_{[0,1]\times\Theta}\pth{r+\lambda^2a}\varphi_m^2-k_m\int_{[0,1]\times\Theta} \varphi_m^2. 
 \end{multline*}
    An integration by parts gives:
  \begin{align*}
    -2\int_{[0,1]\times\Theta} a\varphi_m\dr_x\varphi_m-\int_{[0,1]\times\Theta} \varphi_m^2\dr_xa=0.
  \end{align*}
  We conclude that
  \[\int_{[0,1]\times\Theta}\abs{\dr_x\varphi_m}^2+m\int_{[0,1]\times\Theta}\abs{\nabla_{\theta}\varphi_m}^2\loq \frac{1}{\eta}\pth{(\dabs{r}_{\infty}+\lambda^2\dabs{a}_{\infty})\int_{[0,1]\times\Theta}\varphi_m^2-k_m\int_{[0,1]\times\Theta} \varphi_m^2}.\]
  Since $(k_m)_{m>1}$ is bounded (see Lemma~\ref{lem:estimate_klambda}) and $\dabs{\varphi_m}_{L^2([0,1]\times\Theta)}=1$, the right-hand side of the equation is bounded.
  Thus the family $(\varphi_m)_{m>1}$ is bounded in $W_2^1([0,1]\times\Theta)$.
  Thus, by a compact injection theorem, there is a sequence $m_n\to +\infty$ such that $\varphi_{m_n}$ converges in $L^2([0,1]\times\Theta)$ to a function $\Phi\in L^2([0,1]\times\Theta)$.
  Moreover, the upper bound above implies that as $m\to+\infty$,
  \[\int_{[0,1]\times\Theta}\abs{\nabla_{\theta}\varphi_m}^2\to 0,\]
  so $\Phi=\Phi(x)$ is independent of $\theta$.
  To conclude, we let $\psi\in\mcc^2(\Er)$ be a nonzero principal eigenfunction of the adjoint of the operator defining $\kappa^1(\lambda;a,\overline{r})$, \ie $\psi$ is positive, $1$-periodic and satisfies
    \[{\dr_{x}(a\dr_x\psi)+2\lambda\dr_x(a\psi)+(\overline{r}
      +\lambda^2a-\lambda\dr_xa)\psi}.\]
    We see $\psi$ as function of $\mcc^2(\Er\times\Theta)$ with $\nabla_{\theta}\psi=0$.
    Using $\psi$ as a test function in the weak formulation of Equation \eqref{eq:varphim_mtoinf}, we have:
  \begin{multline*}
    -\int_{[0,1]\times\Theta}a(\dr_x\varphi_m)(\dr_x\psi)-m\int_{[0,1]\times\Theta}\mu{\nabla_{\theta}\varphi_m}\cdot\nabla_{\theta}\psi-2\lambda\int_{[0,1]\times\Theta}a\psi\dr_x\varphi_m\\
    +\int_{[0,1]\times\Theta}(r+\lambda^2a-\lambda\dr_xa)\varphi_m\psi=k_m^{\lambda}\int_{[0,1]\times\Theta}\varphi_m\psi.
  \end{multline*}
  Note that $\nabla_{\theta}\psi=0$. An integration by parts yields 
  \begin{equation*}
    \int_{[0,1]\times\Theta}\varphi_m\dr_x(a\dr_x\psi)+2\lambda\int_{[0,1]\times\Theta}\dr_x(a\psi)\varphi_m
    +\int_{[0,1]\times\Theta}(r+\lambda^2a-\lambda\dr_xa)\varphi_m\psi=k_m^{\lambda}\int_{[0,1]\times\Theta}\varphi_m\psi.
  \end{equation*}
%  Since $\varphi_m$ is bounded in $L^2([0,1]\times\Theta)$, there exists $\Phi\in L^2([0,1]\times\Theta)$ such that up to extraction, $\varphi_m\toweak\Phi$.
  Thus, taking a limit along a subsequence $m_n\to+\infty$ such that $\varphi_{m_n}\to\Phi$ in $L^2([0,1]\times\Theta)$ and $k_{m_n}^{\lambda}$ converges to some limit $k_{\infty}^{\lambda}$, we obtain:
%  \[
%    -\int_{[0,1]\times\Theta}a(\dr_x\Phi)(\dr_x\psi)-2\lambda\int_{[0,1]\times\Theta}a\psi\dr_x\Phi
%    +\int_{[0,1]\times\Theta}(r+\lambda^2a-\lambda\dr_xa)\Phi\psi=k_{\infty}^{\lambda}\int_{[0,1]\times\Theta}\Phi\psi.
%  \]
%  Therefore,  
  \[
    \int_{[0,1]\times\Theta}\Phi\cro{\dr_{x}(a\dr_x\psi)+2\lambda\dr_x(a\psi)
    +(r+\lambda^2a-\lambda\dr_xa)\psi}=k_{\infty}^{\lambda}\int_{[0,1]\times\Theta}\Phi\psi.
    \]
    Since $\Phi$ is independent of $\theta$, the definition of $\psi$ implies:
    \[
    \kappa^1(\lambda;a,\overline{r})\int_{[0,1]\times\Theta}\Phi\psi=k_{\infty}^{\lambda}\int_{[0,1]\times\Theta}\Phi\psi.
    \]
    Since $\varphi_m>0$ for all $m>0$, we must have $\Phi\goq 0$. Since $\dabs{\Phi}_{L^2([0,1]\times\Theta)}=1$, $\Phi$ is not almost everywhere equal to $0$; finally, $\psi>0$.
    Therefore $\kappa^1(\lambda;a,\overline{r})=k_{\infty}^{\lambda}$.
    Thus $\kappa^1(\lambda;a,\overline{r})$ is the only possible limit of a subsequence of $k_m^{\lambda}$ as $m\to+\infty$. Since $(k_m^{\lambda})_{m>0}$ is bounded (Lemma~\ref{lem:estimate_klambda}), we get that $k_m^{\lambda}\to\kappa^1(\lambda;a,\overline{r})$ as $m\to+\infty$.
\end{proof}

\begin{proof}[Proof of the limit $m\to +\infty$ in Theorem~\ref{thm:mutation}.]
  Since $\kappa^1(0;a,\overline{r})>0$, Lemma~\ref{lem:limits_k_muinfty} implies that
  there exist $\delta>0$ and $\underline{m}>0$ such that for all $m>\underline{m}$, we have: $k_m^0>\delta$. 
  Moreover, by Lemma~\ref{lem:estimate_klambda}, for all $m>\underline{m}$ and $\lambda>0$,
  \[k_m^{\lambda}\goq\lambda^2\inf a-\lambda\dabs{\dr_xa}_{\infty}-\dabs{r}_{\infty}.\]
  Therefore,
  \begin{equation}\label{eq:lambda_opti_bdd}
    \lim_{\lambda\to 0}\frac{k_m^{\lambda}}{\lambda}=+\infty,\esp\lim_{\lambda\to +\infty}\frac{k_m^{\lambda}}{\lambda}=+\infty,\qquad\text{uniformly in $m>\underline{m}$}.
  \end{equation}
  By continuity of the function $\lambda\mapsto k_m^{\lambda}$, for all $m>\underline{m}$, there exists $\lambda_m>0$ such that:
  \[c_m=\inf_{\lambda>0}\frac{k_m^{\lambda}}{\lambda}=\frac{k_m^{(\lambda_m)}}{\lambda_m}.\]
  By \eqref{eq:lambda_opti_bdd}, there exist $\lambda_+>\lambda_->0$ such that for all $m>\underline{m}$:
  $\lambda_-\loq\lambda_m\loq\lambda_+$.
  Hence, there exists $\lambda^*>0$ such that up to extraction of a subsequence $m_n\to +\infty$,
  \[\lambda_{m_n}\to \lambda^*\comment{as $n\to+\infty$}.\]
  The functions $\lambda\mapsto k_m^{\lambda}$ ($m>0$) and $\lambda\mapsto\kappa^1(\lambda; a,\overline{r})$ are continuous and convex (see the proof of~\cite[Lemma~3.1]{BHR05-2} for the concavity of $\lambda\mapsto-\kappa^1(\lambda; a,\overline{r})$; the concavity of $\lambda\mapsto -k_m^{\lambda}$ works in the same way).
  Thus 
  \[\lim_{n\to+\infty}\frac{k_{m_n}^{(\lambda_{m_n})}}{\lambda_{m_n}}=\frac{\kappa^1\pth{\lambda^*;a,\overline{r}}}{\lambda^*},\]
  so
  \[\liminf_{n\to+\infty}c_{m_n}\goq c^1(a,\overline{r}).\]
  On the other hand, using the continuity of $\lambda\mapsto\kappa^1(\lambda;a,\overline{r})$, the fact that $\kappa^1(0;a,\overline{r})>0$ and the lower bound in Lemma~\ref{lem:estimate_klambda}, there exists $\overline{\lambda}>0$ such that
  \[c^1(a,\overline{r})=\frac{\kappa^1(\overline{\lambda};a,\overline{r})}{\overline{\lambda}}.\]
  Therefore,
  \[c^1(a,\overline{r})=\lim_{n\to+\infty}\frac{k^{\overline{\lambda}}_{m_n}}{\overline{\lambda}}\goq\limsup_{n\to+\infty}c_{m_n}.\]
  We conclude that $c_{m_n}$ converges and
  \[\lim_{n\to+\infty}c_{m_n}=c^1(a,\overline{r}).\]
  Finally, $(c_{m})_{m>\underline{m}}$ is bounded. Since $c^1(a,\overline{r})$ is the only possible limit of a subsequence $(c_{m_n})$ such that $m_n\to+\infty$, we conclude that the whole family converges to $c^1(a,\overline{r})$.

\end{proof}

\subsection{When the mutation coefficient goes to zero (Theorem~\ref{thm:mutation}, parts~$(ii)$ and~$(iii)$)}\label{ss:m_to_0}

Here, we compute the limit of $c_{m}$ as $m\to 0$ (part~$(ii)$ of Theorem~\ref{thm:mutation}) and we show part~$(iii)$ of Theorem~\ref{thm:mutation}.
%First, we state a lemma analogous to Lemma~\ref{lem:limits_k_muinfty}, but for $m\to 0$ instead of $m\to +\infty$.

\begin{lem}\label{lem:limits_k_muzero}
  Let $\lambda\in\Er$. Under the assumptions of Theorem~\ref{thm:mutation} (part~$(ii)$), as $m\to 0$,
  \[k_m^{\lambda}\to \max_{\theta\in\overline{\Theta}}\kappa^{1}\pth{\lambda;a(\cdot), r(\cdot,\theta)}.\]
\end{lem}

With this lemma, we can conclude the proof of the limiting speed as $m\to 0$. The proof is the same as that of the limit $m\to+\infty$, replacing Lemma~\ref{lem:limits_k_muinfty} by Lemma~\ref{lem:limits_k_muzero}, the limit \enquote{$m\to+\infty$} by the limit \enquote{$m\to 0$}, and $\kappa^1\pth{\lambda;{a},\overline{r}}$ by $\displaystyle\max_{\theta\in\Theta}\kappa^1\pth{\lambda;a,r(\cdot,\theta)}$. 
There only remains to prove Lemma~\ref{lem:limits_k_muzero}.

\begin{proof}[Proof of Lemma~\ref{lem:limits_k_muzero}]
  The principal eigenfunction $\varphi_m^{\lambda}$ solves
\begin{equation*}%\label{eq:eigenvalue_pb_main}
  \left\{
  \begin{aligned}
    \dr_x(a\dr_x{\varphi}^{\lambda}_m)+m\nabla_{\theta}\cdot\pth{\mu\nabla_{\theta}{\varphi}^{\lambda}_m}-2\lambda a\dr_{x}{\varphi}^{\lambda}_m\esp&\\
    +\pth{r+\lambda^2a-\lambda\dr_xa}{\varphi}^{\lambda}_m&=k^\lambda_{m}{\varphi}^{\lambda}_{m},&\text{ in }\Er\times\Theta,\\
    \nu\cdot\nabla_{\theta}{\varphi}_m^{\lambda}&=0,&\text{ over }\Er\times\dr\Theta,\\
    {\varphi}^{\lambda}_m&>0,&\text{ in }\Er\times\Theta,\\
    {\varphi}^{\lambda}_m&\text{ $1$-periodic in $x$}.
  \end{aligned}
  \right.
\end{equation*}
Setting $\eps=m$ and letting $\varphi_{\eps}=\varphi_m^{\lambda}$ brings us to the setting of~\cite[Theorem 1.2]{b23}.
To get back to the the notations there, we replace the slow variable $\theta$ by $y$, the fast variable $x$ by $z$,
the global eigenvalue $k^{\lambda}_m$ by $k_{\eps}$ and the local eigenvalue $k^y$ by $\kappa^1(\lambda;a(\cdot),r(\cdot,\theta))$.
The conclusion is then exactly the result of~\cite[Theorem 1.2]{b23}.
\end{proof}

There remains to prove the last statement of Theorem~\ref{thm:mutation} about the comparison of the two limits $m\to+\infty$ and $m\to 0$.

\begin{proof}[Proof of Theorem~\ref{thm:mutation}, part $(iii)$]
  It is enough to show that for all $\lambda>0$,
  \[
  \kappa^1(\lambda;a,\overline{r})
  \loq
  \sup_{\theta\in\Theta}\kappa^1(\lambda;a(\cdot,\theta),r(\cdot,\theta)).
  \]
  Since $a$ is assumed to be independent of $\theta$, it is equivalent to show that for all $\lambda>0$,
  \[
  \kappa^1(\lambda;a,\overline{r})
  \loq
  \sup_{\theta\in\Theta}\kappa^1(\lambda;a,r(\cdot,\theta)).
  \]
  By~\cite[Proposition 2.4]{Nad10}, the function $g\in L^{\infty}([0,1])\mapsto \kappa^1(\lambda;a,g)$ is convex.
  Hence, by Jensen's inequality, we conclude:
  \begin{align*}
    \kappa^1(\lambda;a,\overline{r})&=\kappa^1\pth{\lambda;a,\frac{1}{\abs{\Theta}}\int_{\Theta}r(x,\theta)\de\theta}\\
    &\loq\frac{1}{\abs{\Theta}}\int_{\Theta}\kappa^1(\lambda;a,r(\cdot,\theta))\de\theta
    \loq\sup_{\theta\in\Theta}\kappa^1(\lambda;a,r(\cdot,\theta)).
  \end{align*}

\end{proof}

\subsection{The optimal mutation coefficient (Theorem~\ref{thm:bio_csq})}\label{ss:proof_bio}

\begin{proof}[Proof of Theorem~\ref{thm:bio_csq}]
  We assume that $a$ is independent of $\theta$ and that $\mu$ is independent of~$x$. 
  We also assume that for some $\theta\in\overline{\Theta}$ and $L>0$, $\kappa^1(0;a,r_{L}(\cdot,\theta))>0$;
  thus by Lemma~\ref{lem:limits_k_muzero} there exist $\delta>0$, $\underline{L}>0$ and $\overline{m}\in\iof{0,+\infty}$ such that for all $m\in\ioo{0,\overline{m}}$,
  we have: $k^0_{m,\underline{L}}>\delta$.
  
  These assumptions ensure that we can apply Theorems~\ref{thm:L_to_inf} and~\ref{thm:mutation} (part~$(ii)$).
  Finally, we assume that \eqref{eq:assumption_hetero} holds. 
  
%  For the simplicity of the proof, we take $\overline{m}=+\infty$, namely, $k^0_{m,\underline{L}}>0$ for all $m>0$. 
%  If on the contrary we cannot take $\overline{m}=+\infty$, then we let
%  \[\overline{m}
%  :=\sup\acc{m \tq \exists \underline{L}>0,\, \exists\delta>0,\,\forall m'\in\ioo{0,m},\, k^0_{m',\underline{L}}>\delta}\]
%  be as large as possible.
%  By Proposition~\ref{ppn:limit_dec_mu}, we have:
%    \[\lim_{m\to \overline{m}\vg m<\overline{m}}\pth{\lim_{L\to+\infty}c_{m,L}}=0,\]
%  and the same proof will work by replacing, in Step 2, $m\to +\infty$ by $m\to \overline{m}$.
%

  \paragraph{Step 1. Proof that $\gamma$ is well-defined and that $\gamma_\mu>0$.} 
  On the one hand, for finite $L>0$, we have by Theorem~\ref{thm:mutation}:
  \[\lim_{m\to 0}c_{m,L}=\max_{\theta\in\overline{\Theta}}c^1(a_L(\cdot),r_L(\cdot,\theta)).\] 
  Thus, using the monotony in $L$ (Proposition~\ref{ppn:inter_L}),
  \begin{equation}\label{eq:lim_Linf_m0}
    \lim_{L\to +\infty}\lim_{m\to 0}c_{m,L}=\max_{\theta\in\overline{\Theta}}\pth{\lim_{L\to+\infty}c^1( a_L(\cdot), r_L(\cdot,\theta))}.
  \end{equation}
  On the other hand, by Proposition~\ref{ppn:L_to_inf},
  \[\lim_{L\to+\infty}c_{m,L}=\pth{\inf_{\lambda>j_m(M)}\frac{j_{m}^{-1}(\lambda)}{\lambda}},\] 
  with
  \[j_{m}(k):=\int_0^1\sqrt{\frac{k-H_m(y)}{a(y)}}\de y\]
  (using the notations $H_m$ and $M$ of Theorem~\ref{thm:L_to_inf}).
  By the Rayleigh formula,
  \[H_{m}(y)=\sup_{\psi\in H^1(\Theta),\,\|\psi\|_{L^2(\Theta)}=1}\pth{\int_{\Theta}r(y,\sigma)\psi^2(\sigma)\de \sigma-m\int_{\Theta}\mu(\sigma)\abs{\nabla\psi}^2(\sigma)\de \sigma}.\]
  Let $m_n\to 0$.
  For fixed $\theta\in\Theta$, we take a sequence $(\psi_n)_n$, $\psi_n\in H^1(\Theta)$ with $\|\psi_n\|_{L^2(\Theta)}=1$, such that $\psi^2_n\toweak \delta_\theta$ in $L^2(\Theta)$ and such that $m_n\dabs{\nabla\psi_n}^2_{\infty}\to 0$. We have for all $y\in\cro{0,1}$ and $n\goq 1$,
  \[H_{m_n}(y)\goq {\int_{\Theta}r(y,\sigma)(\psi_n)^2(\sigma)\de \sigma-m_n\int_{\Theta}\mu(\sigma)\abs{\nabla\psi_n}^2(\sigma)\de \sigma}.\]
  As $n\to+\infty$, the right-hand side of the inequality converges to $r(y,\theta)$.
  This holds for all $\theta\in\Theta$. Hence
  \[\liminf_{n\to+\infty}H_{m_n}(y)\goq\sup_{\theta\in\Theta}r(y,\theta)=\max_{\theta\in\overline{\Theta}}r(y,\theta).\]
  Finally, we also have, for all $y$,
  \[H_{m}(y)\loq \max_{\theta\in\overline{\Theta}}r(y,\theta).\]
  This proves that as $m\to 0$,
  \[H_{m}(y)\to \max_{\theta\in\overline{\Theta}}r(y,\theta).\]
  Therefore, as $m\to 0$, by the dominated convergence theorem,
  \[j_{m}(k)\to j_0(k):=\int_0^1\sqrt{\frac{k-\max\limits_{\theta\in\overline{\Theta}}r(y,\theta)}{a(y)}}\de y.\]
  Using the fact that $\frac{j_m^{-1}(\lambda)}{\lambda}\to +\infty$ uniformly in $m\in(0,1)$ as $\lambda\to 0$ and $\lambda\to +\infty$,
  we get  
  \[\lim_{m\to 0}\lim_{L\to+\infty}c_{m,L}=\inf_{\lambda\goq j_0(M)}\frac{j_0^{-1}(\lambda)}{\lambda};\]
  see e.g. the argument of the proof of Proposition~\ref{ppn:L_to_inf}.
  We obtain (using~\cite[Theorem~2.2]{HNR11}):
  \[\lim_{m\to 0}\lim_{L\to+\infty}c_{m,L}=\lim_{L\to+\infty}c^1\pth{a_L(\cdot),x\mapsto\max_{\theta\in\overline{\Theta}}r_L(x,\theta)}.\]
  By assumption \eqref{eq:assumption_hetero}, we have for all $\theta\in\overline{\Theta}$:
  \[j_0(k)=\int_0^1\sqrt{\frac{k-\max\limits_{\theta\in\overline{\Theta}}r(y,\theta)}{a(y)}}\de y<\int_0^1\sqrt{\frac{k-r(y,\theta)}{a(y)}}\de y.\]
  Therefore, using again~\cite[Theorem~2.2]{HNR11}, we obtain: for all $\theta\in\overline{\Theta}$,
  \[\lim_{L\to+\infty}c^1\pth{a_L(\cdot),x\mapsto\max_{\theta\in\overline{\Theta}}r_L(x,\theta)}> \lim_{L\to+\infty}c^1\pth{a_L(\cdot), r_L(\cdot,\theta)}.\]
  Thus, with \eqref{eq:lim_Linf_m0}, we obtain that the two following limits are well-defined and satisfy:
  \[\lim_{m\to 0}\lim_{L\to+\infty}c_{m,L}> \lim_{L\to+\infty}\lim_{m\to 0}c_{m,L}.\]
  This implies that $\gamma$ is well-defined and satisfies: $\gamma>0$.

  \paragraph{Step 2. There exists an optimal mutation coefficient.}
  Now, we take $\overline{m}$ as large as possible, namely:
  \[\overline{m}:=
  \sup\acc{m'\in\Er\ /\ \exists\underline{L}>0,\,\delta>0,\ \forall m\in(0,m'),\ k^0_{m,\underline{L}}>\delta}
  \in(0,+\infty].\]
  Let $m_1>m_0>0$ with $m_1<\overline{m}$. The family of functions $(m\mapsto c_{m,L})_{L>\underline{L}}$ converges increasingly and pointwise to its limit as $L\to +\infty$ (see Proposition~\ref{ppn:inter_L}), thus the convergence is uniform on $[m_0,m_1]$. Thus for all $\eps>0$, there exists $L_0>\underline{L}$, depending on $m_0$, $m_1$ and $\eps$, such that for all $L_1>L_0$,
  \begin{equation}\label{eq:unif_cv_to_linf}
      \sup_{m\in\cro{m_0,m_1}}\abs{c_{m,L_1}-\lim_{L\to+\infty}c_{m,L}}<\eps.
  \end{equation}
  We set:
  \[\eps:=\frac{1}{2}\pth{\lim_{m\to 0}\lim_{L\to+\infty}c_{m,L}-\lim_{m\to\overline{m},\, m<\overline{m}}\pth{\lim_{L\to+\infty}c_{m,L}}}.\]
  (Recall that $\overline{m}$ may be infinite.)
  By the last point in Proposition~\ref{ppn:limit_dec_mu} and by assumption~\eqref{eq:assumption_hetero} (which implies that there exists $x\in\Er$ such that $\theta\mapsto r(x,\theta)$ is not constant), the function
  \[m\mapsto \lim_{L\to+\infty}c_{m,L}\]
  is strictly decreasing; hence $\eps$ is positive.
Therefore, for $m_0>0$ small enough and $m_1\in(m_0,\overline{m})$ large enough, Equation \eqref{eq:unif_cv_to_linf} implies: for all $m\in (m_1,\overline{m})$, for all $L_1>L_0$,
\[c_{m_0,L_1}>\lim_{L\to+\infty}c_{m,L}\goq c_{m,L_1}.\]
Thus for all $L_1>L_0$,
\[\sup_{m>0}c_{m,L_1}=\sup_{m\in\iof{0,m_1}}c_{m,L_1}.\]
Now, Proposition~\ref{ppn:inter_L} implies that for $L_1> L_0$,
\begin{align*}
  \lim_{m\to 0}c_{m,L_1}\loq \lim_{L\to+\infty}\lim_{m\to 0}c_{m,L}=\lim_{m\to 0}\lim_{L\to+\infty}c_{m,L}-\gamma. 
\end{align*}
Taking $\eps\in\ioo{0,\gamma}$ in \eqref{eq:unif_cv_to_linf} implies:
\[\sup_{m\in\iof{0,m_1}}c_{m,L_1}=\sup_{m\in\cro{m_0,m_1}}c_{m,L_1}.\]
Also, note that $m\mapsto c_{m,L_1}$ is continuous (this can be proved as the continuity of $L\mapsto c_{m,L}$ for fixed $m$, see the first part of the proof of Proposition~\ref{ppn:inter_L}).
Therefore there exists an optimal mutation coefficient $m^*(L_1)\in[m_0,m_1]$. In particular, $m^*(L_1)>0$.

\paragraph{Step 3. Lower bound for the maximal speed.}
Let $\eps>0$.
By virtue of Proposition~\ref{ppn:limit_dec_mu}, there exist $m_1>0$ so small that for all $m\in\iof{0,m_1}$,
\[\lim_{L\to+\infty}c_{m,L}\goq \lim_{m\to 0}\lim_{L\to+\infty}c_{m,L}-\eps/2.\]
By \eqref{eq:unif_cv_to_linf} with $\eps$ replaced by $\eps/2$, for all $L_1\goq L_0$, for all $m\in\iof{0,m_1}$,
\begin{align*}
  c_{m,L_1}&\goq \lim_{L\to+\infty}c_{m,L}-\eps/2\goq\lim_{m\to 0}\lim_{L\to+\infty}c_{m,L}-\eps\\
  &\goq\lim_{L\to+\infty} \lim_{m\to 0}c_{m,L}+\gamma-\eps\\
  &\goq \lim_{m\to 0}c_{m,L_1}+\gamma-\eps.
\end{align*}
This concludes.

\end{proof}

\subsection{Other properties (Propositions~\ref{ppn:limit_dec_mu} and~\ref{ppn:persistence_decr})}\label{ss:other}

\begin{proof}[Proof of Proposition~\ref{ppn:limit_dec_mu}]
  We first prove that
  \[m\mapsto \lim_{L\to+\infty}c_{m,L}\]
  is nonincreasing (the limit exists thanks to the assumptions and Theorem~\ref{thm:L_to_inf}).
  By Proposition~\ref{ppn:L_to_inf}, it is enough to show that for all $k\goq M$, $m\mapsto j_m(k)$ is nondecreasing, which is implied by the property that for all $x\in\Er$, $m\mapsto H_m(x)$ is nonincreasing . 
By the Rayleigh formula, we have:
  \[H_m(x)=\max_{\phi\in H^1(\Theta),\,\|\phi\|_{L^2(\Theta)}=1}\pth{-m\int_{\Theta}\mu\abs{\nabla\phi}^2+\int_{\Theta}r(x,\cdot)\phi^2}.\]
  Take $m_1,m_2\in I$ with $m_1<m_2$, and let $\phi_2\in H^1(\Theta)$ with $\phi_2>0$ and $\|\phi_2\|_{L^2(\Theta)}=1$ be the principal eigenfunction associated to the principal eigenvalue $H_{m_2}(x)$. Then:
  \begin{align*}
    H_{m_1}(x)&\goq{-m_1\int_{\Theta}\mu\abs{\nabla\phi_2}^2+\int_{\Theta}r(x,\cdot)\phi_2^2}\\
    &\goq{-m_2\int_{\Theta}\mu\abs{\nabla\phi_2}^2+\int_{\Theta}r(x,\cdot)\phi_2^2}=H_{m_2}(x).
  \end{align*}
  Thus for all $x\in\Er$, the function $m\mapsto H_m(x)$ is nonincreasing. Moreover, if $\phi_2$ is not constant, the second inequality is strict, so the function $m\mapsto H_m(x)$ is strictly decreasing. Hence the function
  \[m\mapsto\lim_{L\to+\infty}c_{m,L}\]
is nonincreasing, and strictly decreasing if there exists $x\in\Er$ such that $\theta\mapsto r(x,\theta)$ is not constant. 

Now, let us prove that
  \[m\mapsto \lim_{L\to 0}c_{m,L}\]
  is nonincreasing (the limit exists thanks to the assumptions and Theorem~\ref{thm:L_to_0}). This comes from the variational formulation in Theorem~\ref{thm:L_to_0},
  which says:
  \[\lim_{L\to 0}c_{m,L}=\sup_{\psi\in H^1(\Theta)\vg\dabs{\psi}_{L^2(\Theta)}=1,\, R_{\psi}>0}2\sqrt{A_{\psi}R_{\psi}},\]
  where, for all $\psi\in H^1(\Theta)$, $A_{\psi}$ is independent of $m$ and
  \[R_{\psi}:=\int_{\Theta}\pth{\widetilde{r}^A\psi^2-m\widetilde{\mu}^A\abs{\nabla\psi}^2}\]
  is nonincreasing in $m$. We conclude in the same way as above.

\end{proof}

\begin{proof}[Proof of Proposition~\ref{ppn:persistence_decr}]
  By the Rayleigh formula,
  \[k^0_m=\sup_{\phi}\pth{-\int_{[0,1]\times\Theta}a\abs{\nabla_x\phi}^2-m\int_{[0,1]\times\Theta}\mu\abs{\nabla_{\theta}\phi}^2+\int_{[0,1]\times\Theta}r\phi^2},\]
  where the $\sup$ is taken on the $\phi\in H^1_{loc}(\Er\times\overline{\Theta})$ such that $\|\phi\|_{L^2([0,1]\times\Theta)}=1$ and $\phi$ is $1$-periodic in $x$.
  From this expression, we get that $m\mapsto k^0_{m}$ is nonincreasing.
\end{proof}

\section*{Acknowledgements}
The author thanks Raphaël Forien, François Hamel and Lionel Roques for their advice and support.
The author was supported by the ANR project ReaCh ({ANR-23-CE40-0023-01}) and by the Chaire Modélisation Mathématique et Biodiversité (École Polytechnique, Muséum national d’Histoire naturelle, Fondation de l’École Polytechnique, VEOLIA Environnement).

\printbibliography

\end{document}